\documentclass[a4paper,12pt]{amsart}
\usepackage[utf8]{inputenc}
\usepackage{amsthm}
\usepackage{amsmath}
\usepackage{bm}
\usepackage{enumitem}
\usepackage{amsfonts}
\usepackage{amssymb}
\usepackage{mathtools}
\usepackage{appendix}
\usepackage{mathrsfs}
\usepackage{setspace}
\usepackage{comment}
\usepackage{xcolor}
\usepackage{todonotes}
 \usepackage{ulem}
\usepackage{thmtools} 
\usepackage[foot]{amsaddr} 
\usepackage{appendix}
\RequirePackage[a4paper,top=2.54cm,bottom=2.54cm,left=1.90cm,right=1.90cm,%
                headsep=1em,includehead,includefoot]{geometry}
\usepackage[pdfdisplaydoctitle,colorlinks,breaklinks,urlcolor=blue,linkcolor=blue,citecolor=blue]{hyperref} 
\usepackage[nameinlink,capitalise]{cleveref}
\usepackage{soul}

\providecommand{\U}[1]{\protect\rule{.1in}{.1in}}
\newtheorem{theorem}{Theorem}

\newtheorem{corollary}[theorem]{Corollary}

\newtheorem{assumption}[theorem]{Assumption}

\newtheorem{lemma}[theorem]{Lemma}

\newtheorem{remark}[theorem]{Remark}

\newcommand{\R}{\mathbb{R}}
\newcommand{\Div}{\text{div}}
\title[Hasegawa-Mima equation from compressible Euler-Boltzmann system]
{Hasegawa-Mima equation as limit of a magnetized pressureless Euler-Boltzmann system}

\author[Federico Butori]{Federico Butori}
\address{Scuola Normale Superiore, Piazza dei Cavalieri, 7, 56126 Pisa, Italy}
\email{\href{mailto:federico.butori at sns.it}{federico.butori at sns.it}} 
\author[Franco Flandoli]{Franco Flandoli}
\address{Scuola Normale Superiore, Piazza dei Cavalieri, 7, 56126 Pisa, Italy}
\email{\href{mailto:franco.flandoli at sns.it}{franco.flandoli at sns.it}}
\author[Yassine Tahraoui]{Yassine Tahraoui}
\address{Scuola Normale Superiore, Piazza dei Cavalieri, 7, 56126 Pisa, Italy}
\email{\href{mailto:yassine.tahraoui at sns.it}{yassine.tahraoui at sns.it}}
\date\today

\begin{document}
	\begin{abstract} 
We prove that the Hasegawa-Mima equation (HME) is the limit of a magnetized pressureless Euler-Boltzmann system in a special scaling regime. The proof relies on the derivation of a uniform estimate independent of the scaling parameter, using the particular structure of the singular system and the identification of a combination of the unknowns that compensates for fast oscillations in time. We also show that HME allow us to investigate a conditional property of inverse energy cascade and drift waves.
     \end{abstract}
	\maketitle
\tableofcontents
\section{Introduction}
The Hasegawa-Mima equation  \eqref{HM-eqn-Intro} is a reduced plasma model of great value
for the understanding of turbulence in plasma, see \textit{e.g.} \cite{HasMim,HorHas}.
It describes, up to an appropriate normalization, the fluctuations of the electrostatic potential $\phi$ on the poloidal plane, under the action of an average (logarithmic) density gradient $\nabla \log n_0$ and under the assumption of a strong  guide field directed orthogonally to the poloidal plane. 
 The equation, in the unknown $\phi:\left[  0,T\right] \times \R^2  \rightarrow
\mathbb{R}$, reads in its dimensionless form
\begin{align}\label{HM-eqn-Intro}
\tag{HME}
	\begin{cases}
		&\partial_t(\phi-\Delta \phi)-\nabla^\perp \phi\cdot\nabla\left( \log n_0- \Delta \phi\right) =0 \quad  \text{ in } [0,T]\times \R^2, \\
		&\phi_{\vert t=0}=\phi_0, 
	\end{cases}
\end{align}
where  $n_{0}$ is a positive function, independent of time  representing a background density,  $\phi_0$ is the initial datum, $\nabla^\perp \phi=(\partial_2\phi,-\partial_1\phi)$  and  the final time horizon $T>0$.
Beside a wide
literature in physics, a number of works exist dealing with rigorous
properties of the \eqref{HM-eqn-Intro}, see for instance \cite{GuoHan,butori2026ergodicity,CaoFar,Flandoli-Tahraoui}. No one of these works treat the derivation of \eqref{HM-eqn-Intro} from more
fundamental plasma equations and  the rigorous derivation is the purpose of the
present paper, see \autoref{Theorem-main}. 
In this work, we will prove that \eqref{HM-eqn-Intro} describes the singular limit $\varepsilon \rightarrow 0$ of the pressureless Euler-Boltzmann system 
\begin{equation*}\label{eul-bol_intro}\tag{$\varepsilon$-EB}
    \begin{cases}
        \partial_t \rho_\varepsilon +\Div(\rho_\varepsilon u_\varepsilon)=0\\
        \partial_t(\rho_\varepsilon u_\varepsilon) + \Div(\rho_\varepsilon u_\varepsilon\otimes u_\varepsilon)=\frac1 \varepsilon \rho_\varepsilon\left(-\nabla \phi_\varepsilon + u_\varepsilon^\perp\right) \\
        \rho_\varepsilon(x)= n_0^\varepsilon e^{\varepsilon \phi_\varepsilon}. \\
        {\rho_{\varepsilon}}_{|t=0} =\rho_{\varepsilon, 0}= n_0^\varepsilon e ^{\varepsilon \tilde \phi_0}, \  {u_{\varepsilon}}_{|t=0} =u_{0}
    \end{cases}
\end{equation*}
where $\nabla \log n_0:=\lim_{\varepsilon\rightarrow0}\varepsilon^{-1} \nabla \log n_0^\varepsilon$. In particular, we will find along the way that, in an appropriate sense, $\rho_\varepsilon$ converges to a constant density profile (c.f. \autoref{rem:density})
and $\lim_{\varepsilon\rightarrow 0} u_\varepsilon=- \nabla^\perp \phi$, so that the system is in fact fully describes only by the electrostatic potential $\phi$. Finally, as a consequence of the singular term in the momentum equation, the initial condition $\phi_0$ of \eqref{HM-eqn-Intro} will be given by 
$(I-\Delta)^{-1}(\tilde \phi_{0}- \nabla^\perp \cdot u_{0})$, where $\tilde \phi_0$ is the given initial electric potential  for \eqref{eul-bol_intro} thus it coincides with $\tilde \phi_0$ only in the case of well prepared data $u_0=\nabla^\perp \tilde \phi_0$. 
The derivation of \eqref{HM-eqn-Intro}, on a physical ground, starting from even more fundamental systems, is reported in many works as those
quoted above. It is based on a few initial assumptions: cold-ions, 
adiabatic electrons, quasi-neutrality and the two-dimensional reduction, which simplify the starting plasma model, plus an ordering condition, on the relevant scales of the system (the parameter $\varepsilon$).  However, no rigorous result has been proved yet, to the best of the authors' knowledge.
The heuristic mechanism of proof by means of the order parameter $\varepsilon $
is not so trivial, compared to more intuitive scaling limits like the one of the strong magnetic field regime. In  \autoref{sec:heuristic} below,
we illustrate a proper mathematical reformulation of the heuristic proof,
before we proceed to the rigorous analysis in the rest of the paper.
Reducing complex plasma systems to equations like the \eqref{HM-eqn-Intro} one is very
important, because the simplified models allow us to investigate questions
that are too difficult or obscure for the fluid dynamic system \eqref{eul-bol_intro}. In \autoref{sec:cascade}, we
show two examples about  a conditional property of
inverse energy cascade and drift waves. \\

In the rest of this section we give an overview of the main ideas and challenges to establish the limit rigorously and we collect some physical considerations that motivate the choice of system \autoref{eul-bol_intro} as a starting point. 

 \subsection{The strategy of the proof}\label{sec:heuristic}
 The system \eqref{eul-bol_intro} can be formally recast into a system of two equations in two unknowns. Since $\partial_t \phi_\varepsilon  = \varepsilon^{-1}\partial _t \rho_\varepsilon  /\rho_\varepsilon $, assuming $\rho_\varepsilon  >0$, we can arrive with simple computations at the following equations for $(\phi_\varepsilon , u_\varepsilon )$: 
        \begin{align}\label{Euler-rho-phi-Main}
        \tag{$\varepsilon$-PV}
	\begin{cases}
	&\partial_{t}\phi_\varepsilon + \frac{1}{\varepsilon}\operatorname{div}    u_\varepsilon +u_\varepsilon \cdot \left(\nabla \phi_\varepsilon  + g_\varepsilon \right)  =0, \quad 0<t<T, \quad x\in \mathbb{R}^2,\\[0.0cm]
&\partial_{t}  u_\varepsilon  + u_\varepsilon \cdot \nabla u _\varepsilon   
=\frac{1}{\varepsilon}\left(-\nabla\phi_\varepsilon 
  +  u_\varepsilon ^{\perp}\right).
	\end{cases}
	\end{align}
where we have introduced the auxiliary function 
$$g_\varepsilon(x) := \frac{1}{\varepsilon}\nabla \log n_0^\varepsilon$$
Assume for the moment that $(\phi_{\varepsilon},u_{\varepsilon})_\varepsilon$ are bounded  in suitable sense (see \autoref{Sec:unif-estimates} and \autoref{Prop-Hs-estimate}) and thus that they converge (up to subsequences) in a suitable weak sense to some pair $\left(
\phi,u\right) $. From the momentum equation, we may argue that the singular penalization yields $$-\varepsilon\partial_{t}u_{\varepsilon}^{\perp} -\varepsilon\left(  u_{\varepsilon}  \cdot\nabla
u_{\varepsilon} \right)  ^{\perp}\to 0$$  and thus
\begin{equation}\label{limit-phi-u}
u^\perp   -\nabla\phi =0 ,
\end{equation}
and since $(u^\perp)^\perp=-u$, we obtain $u=- \nabla^\perp \phi$. Next, assuming $g_\varepsilon \rightarrow g:= \nabla G$ in an appropriate sense and 
passing to the limit in the first equation of \eqref{Euler-rho-phi-Main}, yields
\[
 \partial_{t}\phi +u  \cdot\nabla \phi =\operatorname{div}\left(  \partial_{t}u^{\perp}  +  u  \cdot\nabla u^\perp 
\right)  -u \cdot g  .
\]
Finally, replacing $u=-\nabla^{\perp}\phi$,
we get the Hasegawa-Mima equation: 
\begin{align*}
    \partial_t(\phi-\Delta \phi)+\nabla^\perp \phi\cdot\nabla(\Delta \phi)-\nabla^\perp \phi\cdot g=0.
\end{align*}\\

\paragraph{\textit{Rigorous proof and difficulties.}}
The naive argument of above contains several crucial difficulties. The first one is that
the well posedness of \eqref{HM-eqn-Intro} requires a good degree of regularity, which
is not encoded in the conserved quantities of the Euler-Boltzmann system we would like to start with. Even the definition of $\phi $
in terms of ion density $\rho $ by Maxwell-Boltzmann formula requires
non-trivial properties like $\rho >0$ (not only $\geq 0$). 
 Moreover, the heuristic derivation of the Hasegawa--Mima equation \eqref{HM-eqn-Intro} from \eqref{eul-bol_intro} suggests working with smooth solutions rather than merely weak solutions. Indeed, the  analysis relies on strong compactness and on the possibility of performing nonlinear manipulations involving derivatives of the unknowns. In particular, one needs to justify the passage to the limit in nonlinear expressions involving products of derivatives. 
We overcome this difficulty by working with smooth solutions of the Euler system, provided by
the Kato theory of symmetric hyperbolic systems \cite{kato1975cauchy}. The systems \eqref{Euler-rho-phi-Main} and  \eqref{eul-bol_intro} are a first-order quasilinear hyperbolic systems and
 for sufficiently smooth initial data $H^s(\R^2),s>2$, the classical theory of symmetric hyperbolic systems (see  \textit{e.g.} \cite{kato1975cauchy}) yields the local-in-time existence and uniqueness of  $H^s(\R^2)$-solution. This, however, restricts the result on a local time interval.\\
 
At the same time, though, this theory does not provide a lifespan that is uniform with respect to the singular parameter $\varepsilon$. Indeed, a direct energy estimate involves coefficients of order $1/\varepsilon$. As a consequence, a standard Grönwall argument leads, in general, to a local existence time $T_\varepsilon$ that depends on $\varepsilon$ and may shrink to zero as $\varepsilon\to0$. Obtaining a lifespan that is uniform in $\varepsilon$ therefore requires exploiting the specific structure of \eqref{eul-bol_intro} and identifying suitable cancellations in the singular terms, rather than relying solely on the general hyperbolic theory. This can be considered the first major difficulty to be solved, and the aim of \autoref{Sec:unif-estimates} is to prove that the lifetime can be chosen independently of $\varepsilon$, see \autoref{Prop-Hs-estimate}.\\

Again inspired by special properties of the limit \eqref{HM-eqn-Intro}, we introduce an auxiliary quantity, a combinations of the unknowns of \eqref{eul-bol_intro}, whose time derivatives remain uniformly controlled due to suitable cancellations.  As a consequence, we obtain strong convergence of an appropriate quantity, which is sufficient to get the limit \eqref{HM-eqn-Intro}.
 Let us explain this in more details: consider the quantity $q_\varepsilon=
\phi_\varepsilon+ \nabla^\perp \cdot u_\varepsilon$ (where $\nabla^\perp = (\partial_2, - \partial_1)$). By differentiating in time, we can see that $q_\varepsilon$ solves the following system 
\begin{align*}
    \begin{cases}
    &\partial_t q_\varepsilon
+
u_\varepsilon\cdot\nabla q_\varepsilon
+(\nabla^\perp \cdot u_\varepsilon)
\operatorname{div}u_\varepsilon
+
u_\varepsilon\cdot g_\varepsilon
=0,\\[0.15cm]
&q_\varepsilon=
\phi_\varepsilon+\nabla^\perp\cdot u_\varepsilon, \\
&\big(\phi_\varepsilon, u_\varepsilon\big)    \big|_{t=0} = \big(\phi_0, u_{0}\big).
\end{cases}
\end{align*}
Crucially, all the singular terms have disappeared from the equation, thus it is possible to prove $\varepsilon$-independent estimates for $q_\varepsilon$. Such estimates provide the temporal compactness required to justify rigorously the convergence towards the Hasegawa--Mima equation \eqref{HM-eqn-Intro}. Indeed, up to subsequences, there exists a limit $q$ such that $q_\varepsilon \rightarrow q$ in some topology (see \autoref{strong-cv-lemma}). Combining this with the identification of the limits for $\phi_\varepsilon$ and $u_\varepsilon$ (recall $u = - \nabla^\perp \phi$ from \eqref{limit-phi-u}) yields
$$q = \phi + \nabla^\perp \cdot u= \phi - \Delta \phi.$$
Moreover, the strong convergence of $q_\varepsilon$ allows to pass to the limit in the nonlinear term $u_\varepsilon\cdot \nabla q_\varepsilon$. For the defect term $(\nabla^\perp \cdot u_\varepsilon)
\operatorname{div}u_\varepsilon$, altough $\Div (u_\varepsilon)$ weakly converges to zero (see \autoref{lemma-0limit}), this does not suffice for controlling the convergence of the product, which however, thanks to the specific structure of the equation
satisfies, in distribution,
\begin{equation*}
    (\nabla^\perp \cdot u_\varepsilon)\operatorname{div}(u_\varepsilon) = \frac{\varepsilon}{2}\left[ \partial_t|\nabla^\perp \cdot u_\varepsilon|^2 + \Div(u_\varepsilon |\nabla^\perp \cdot u_\varepsilon|^2) + \Div (u_\varepsilon) |\nabla^\perp \cdot u_\varepsilon|^2\right]
\end{equation*}
and thus weakly converges to zero as $\varepsilon\rightarrow0$, since all terms in the square bracket are uniformly controlled, thanks to our a-priori estimates. One then readily obtains equation \eqref{HM-eqn-Intro}.
Thus, the key point is not to control $\partial_t\phi_\varepsilon$ and $\partial_tu_\varepsilon$ separately, but rather the time derivative of $q_\varepsilon$, which encodes the slow dynamics of the system, see \autoref{Sec:derivtionHM} for more details.
\subsection{The low amplitude--long wavelength regime}\label{sec: physics}
In many situations where a strong and constant guiding magnetic field is applied to a three dimensional plasma, a good description of the dynamics of the plasma in the poloidal plane (that orthogonal to the direction of the magnetic field) can be given by two dimensional fluid models.
One that is particularly relevant for regimes of nuclear fusion is the following 2D pressureless Euler-Boltzmann system:
\begin{align}\label{Euler-new-Main}
\tag{EUL-BOL}
	\begin{cases}
	&\partial_{t}n+\operatorname{div}\left(  U
n\right)     =0,\\[0.0cm]
&\partial_{t}\left(  Un\right)  +\operatorname{div}%
\left(  n U\otimes U\right)   
=\frac{\mathrm{q}}{m_{i}}n\left(
-\nabla  \Phi  +B U^{\perp}\right),\\[0.0cm]
&n\left(
x,t\right)  =\tilde N_{0}\left(  x\right)  e^{\frac{\mathrm{q}}{T_{e} }%
\Phi\left(  x,t\right)  }
	\end{cases}
	\end{align} 
    where $n$ denotes the plasma density, $U$ the plasma velocity, $\Phi$ the electrostatic potential, $\mathrm{q}, m_i$ are the charge and the mass of ions, $T_e$ is the temperature of electrons, $\mathrm{q}$ is the (positive) electric charge (while $\mathrm{q}_e=-\mathrm{q}$ is the negative electron charge), $B$ is the intensity of the guide magnetic field and we denoted $U^\perp = (U_2 , - U_1)$. 
    Here we consider a model for ions in a fixed background of (massless) electrons in a quasi-neutral regime (the more physics oriented reader could be surprised by the absence of a Poisson term in our Maxwell-Boltzmann. This is indeed the \textit{quasineutral} regime. See also \autoref{remark-cold-ions} for a broader discussion). Consequently, $\tilde N_0$ represent a background density and we assume that electrons are \textit{thermalized}, i.e. they have reached their thermodynamic equilibrium, expressed by the third equation, the Maxwell-Boltmann law.
    Finally we supplement the system with initial conditions $\rho\left(  0,x\right)     =\tilde N_{0}\left(  x\right)
e^{\frac{\mathrm{q}}{T_{e}}{\Phi}_{0}\left(  x\right)  }$,  $ U\left(  0,x\right)   ={U}_{0}\left(
x\right).$ 
A key feature of this model is the absence of pressure in the momentum equation. This is an assumption on the system known in physics as `cold ions' or `cold plasma' (see \cite[Ch 4]{chen2018introduction}).  We refer to \autoref{remark-cold-ions} for more details about the derivation of the Euler-Boltzmann model \eqref{Euler-new-Main} from more fundamental ones. 
 
\begin{remark}
Here, as in most of the physics literature, we assume
$T_{e}\left(  x\right)  $ is constant in space and time and equal to $T_{e}$. A
rigorous derivation with space-dependent temperature would be an interesting extension of our main theorem, but outside the scope of the present manuscript. The final result has been indicated and
studied for instance by \cite{Nyc}.
\end{remark}
    In order to introduce the special regime in which the Hasegawa-Mima equation \eqref{HM-eqn-Intro} are obtained from \eqref{Euler-new-Main}, let $L$, $\tau$, $\hat u$, $\hat \phi$, $\hat n$, denote reference length, time, velocity, electrostatic and density. We introduce the dimensionless variables
\begin{equation}
\label{scaling}
x'=\frac xL,
\qquad
t'=\frac{t}{\tau} ,
\qquad
u=\frac {U}{\hat u} , \qquad \phi= \frac{\Phi}{\hat \phi}, \qquad \rho=\frac{n}{\hat n}, \qquad \tilde n_0=\frac{\tilde N_0}{\hat n}.
\end{equation}
We also assume that the background density $\tilde n_0$ oscillates at a characteristic wavelength $\ell$, in particular we impose
\begin{equation}\label{ass:scaling_rho}\tag{A1}
    \nabla \log \tilde n_0= \ell^{-1} \nabla \log n_0
\end{equation}

The regime of interest can then be expressed introducing an ordering parameter $\varepsilon>0$ and assuming that
\begin{equation}\label{ass:scaling_eps} \tag{A2}
    \varepsilon = \frac{\mathrm{q}\hat \phi }{T_e} = \frac{m_i}{\mathrm{q}B\tau} =  \frac{m \hat uL}{\mathrm{q}\tau \hat{\phi}} = \frac{L}{\ell}, \qquad L= \tau \hat u = \frac{(T_em_i)^{1/2}}{\mathrm{q}B}.
\end{equation}
This assumptions corresponds to a regime of small amplitude and long wavelength plasma oscillations, we give more details in the next paragraph. 
Rescaling \eqref{Euler-new-Main} with this choice of dimensionless constants gives \eqref{eul-bol_intro}.
 The main contribution of the present work is to show that equation \eqref{HM-eqn-Intro} arises in the limit $\varepsilon\rightarrow 0$ from the system \eqref{eul-bol_intro}. Consequently, our result justifies that in the small $\varepsilon$ regime, the electrostatic potential $\Phi_\varepsilon$ can be described by $\varepsilon \frac{T_e}{q}\phi$, where $\phi$ is given by a solution of \eqref{HM-eqn-Intro}. 

 \begin{remark}
     The meaning of the ordering relations \eqref{ass:scaling_rho} and \eqref{ass:scaling_eps} might be confusing for a more mathematic oriented reader. They can be recast as assumptions on $n_0$ and smallness of the initial data and we will explain how to do it later in \autoref{rem:scaling}. Under those assumptions, the link between \eqref{Euler-new-Main} and \eqref{eul-bol_intro} is obtained just by rescaling. 
 \end{remark}
The parameter $\varepsilon$ prescribes precise relations between the different reference scales of the plasma. The cyclotronic frequency of the plasma is given by $\omega_{c}=\mathrm{q}B/m$. The first relation of \eqref{ass:scaling_eps} gives $\omega_{c}\tau = \varepsilon^{-1}$ which means a reference time scale much larger than the cyclotronic time scale. Next the reference electrostatic energy is given by $\mathcal{E}_{el} = \mathrm{q}\hat \phi$, while the reference Kinetic energy by $m \hat u^2$ and the thermal energy (or temperature) of the electrons is given by $T_e$, thus the second relation in \eqref{ass:scaling_eps} gives the smallness of
the ratio $\mathcal{E}_{el}/T_e$ and finally the third relation in \eqref{ass:scaling_eps} (plus $L= \tau \hat u$) gives $\mathcal{E}_{kin}/\mathcal{E}_{el} = \varepsilon$. This is the reason why this regimes takes the name of \textit{low amplitude}. Finally, the relationship $L/\ell\sim \varepsilon$ corresponds to \textit{long wavelength} oscillations of the background density, namely on a scale much larger than the reference scale $L$.
Finally, regarding spatial scales, let us introduce the \textit{ion sound Larmor radius} (often called the ion sound gyroradius) $\rho _{s}$ which is the characteristic length scale of an ion moving at the plasma sound speed in a magnetic field. It is given by $\rho _{s}:= \sqrt{mT_e} /\mathrm{q}B$ and we see that the last relation of \eqref{ass:scaling_eps} express exactly $L/\rho_s=1$. This ratio is the only dimensionless quantity that is not small with $\varepsilon$.
Thanks to these remarks, we can draw a direct equivalence between the ordering relation \eqref{ass:scaling_eps} and the ordering found in \cite[Equation 2.1]{HH} (notice that we don't have the second relation here because we assumed that the plasma is two dimensional from the beginning). \\

Let us end this section by reintroducing the physical constants in the Hasegawa-Mima equation. By performing the inverse transformations of $\eqref{scaling}$ on \eqref{HM-eqn-Intro}, (recalling that $L=\rho_s$ and $\tau^{-1}= \varepsilon \omega_c$), that is defining 
$\Phi(t,x)= \varepsilon\frac{T_e}{q}\phi( \varepsilon \omega_ct, \frac{x}{\rho_s})$ and a function $N_0(x)$ such that $\nabla \log N_0(x) = \varepsilon \nabla \log n_0(\frac{x}{\rho_s})$ we get the equation in dimensional variables
\begin{align}\label{HM-D}
\tag{HME-dim}
	\begin{cases}
		&\frac{1}{\omega_{c}}\partial_t(\frac{\mathrm{q}\Phi}{T_e}-\rho_s^2\Delta \frac{\mathrm{q}\Phi}{T_e})+\rho^4_s\nabla^\perp \frac{\mathrm{q}\Phi}{T_e}\cdot\nabla\Delta \frac{\mathrm{q}\Phi}{T_e} - \rho^2_s\nabla^\perp \frac{\mathrm{q}\Phi}{T_e}\cdot \nabla \log N_0=0 \quad  \\
		&\Phi_{\vert t=0}=\varepsilon\frac{T_e}{q}\phi_0( \frac{x}{\rho_s}).
	\end{cases}
\end{align}
 Given our derivation of the previous sections, the reader should not be surprised that, in dimensional variables, the (drift) velocity of the plasma is given by $U=-\frac{1}{B}\nabla^\perp \Phi$ (the so-called $E\times B$ \textit{drift}).
 Equation \eqref{HM-D} contains many physical quantities. Let us recall that 
 $$\omega_c = \frac{\mathrm{q}B}{m}, \qquad \rho _{s}:= \frac{\sqrt{mT_e} }{\mathrm{q}B}$$
 From an engineering point of view, ion mass $m$ and charge $\mathrm{q}$ are given, but in context of fusion plasma, also the temperature $T_e$ is fixed (roughly $10^8$ degrees, the minimum value to allow for fusion to take place). In this regard, the only free engineering parameter is the intensity of the magnetic field $B$. It can be useful for intuition to rewrite \eqref{HM-D} underlying only the dependence on $B$ in the various constants, i.e. assuming for simplicity $\mathrm{q}=m=T_e=1$. The equation that we get is then 
 \begin{align}\label{HM-B}
\tag{HME-B}
		&\partial_t(\Phi-\frac{1}{B^2}\Delta \Phi)+\frac{1}{B^3}\nabla^\perp \Phi\cdot\nabla\Delta \Phi - \frac{1}{B}\nabla^\perp \Phi\cdot \nabla \log N_0=0 
\end{align}

\subsection{Link with Vlasov-Poisson equations}\label{remark-cold-ions}
The system \eqref{Euler-new-Main} is usually derived from a kinetic description of the plasma, see for instance \cite{Bra}. The most complete description would requires two sets of coupled equations, one for the electrons density and one for the ions' (see the seminal work \cite{bouchut1991global}). In the regime of small electron to ions mass ratio, electrons tend to \textit{thermalize}, i.e. their distribution tends to a Maxwellian, indepepndent of the ions' distribution. This allows to reduce the model to a single equation for ions with a Maxwell-Boltzmann law for the background of electrons (see \cite{grenier2020derivation} for a rigorous proof for a two-fluids model and \cite{BARDOS201894} for a kinetic model). Consequently, the electrostatic potential should follow a Poisson law, which reads, in the case of thermalized or massless electrons 
$$\varepsilon_0\Delta \phi = \tilde N_0 e^\frac{q\phi}{T_e} - n.$$
When passing to the dimensionless variables employing the scaling \eqref{scaling}, one obtains the relation 
$$\delta^2\Delta \phi = n_0 e^\phi - \rho.$$
where the parameter $\delta^2 = \varepsilon_0T_e/(q\hat nL^2) $ is the square of the ratio between the Debye length and the observation length $L$. In laboratory plasmas this ratio is usually extremely small, hence it is reasonable to assume \textit{quasineutrality}, i.e.  $\delta=0$ and recover the Bolzmann law from \eqref{Euler-new-Main}. 
In general one could add to relations \eqref{ass:scaling_eps}, the ordering $\delta^2= \varepsilon^2$ and treat also the quasineutral limit as a higher order correction. 
Since we believe that this generalization would not add much insight into the derivation of the Hasegawa Mima equations, at the cost of complicating many of the arguments employed, we have decided to defer it to a future work. 
All in all, the distribution $f(t,x,v)$ of ions in a background of thermalized electrons, in the quasineutral regime evolves accordingly to the Vlasov-Boltzmann system
\begin{align}\label{eq:vlasov}
\begin{cases}
    \partial_t f + v\cdot \nabla_x f + \frac{q_i}{m_i}(- \nabla_x \phi + Bv^\perp)\cdot \nabla_v f =0 \\
    n(t,x)= \int f(t,x,v) dv, \qquad 
    n(t,x) = n_0(x) e^{\frac{q}{T_e}\phi(t,x)}.
    \end{cases}
\end{align}
The well posedness of the above system (including also the Poisson term in the equation for the electrostatic potential) was recently discussed in \cite{griffin2020global}, while the \textit{quasineutral limit} leading to \eqref{eq:vlasov} was established in \cite{griffin-pickering_singular_2020a}. The link between \eqref{eq:vlasov} and \eqref{Euler-new-Main} is given by the \textit{monokinetic assumption}, i.e. if $(n, U)$ is a solution of \eqref{Euler-new-Main}, then $f(t,x,v):=n(t,x) \delta_{v=U(t,x)}$ is a measure-valued solution of \eqref{eq:vlasov}, and \textit{vice versa}; in a more abstract framework, this analogy can even be taken further as discussed in \cite{brenier_dual_1993,brenier_least_1989,zakharov_benney_1980}. In particular, the special monokinetic structure is conserved in the evolution and yields the vanishing of the ion's pressure. On the other hand, measure-valued solutions of the Vlasov system \eqref{eq:vlasov} have been considered several times in the literature, often describing fluid system with current/vortex sheets, or in some cases, representing well prepared data which allow to avoid certain instabilities present in kinetic models, see \cite{baradat_nonlinear_2020, brenier_homogenized_1997, grenier_oscillations_1996}. 
Note that solution with Dirac delta velocities ideally correspond to Gaussian velocity profiles with zero temperature, from which the denomination of `cold ions'. 

\subsection{Related models and singular limits}
There are very important examples of rigorous derivation of various other
reduced fluid plasma models, mostly based on the passage from the kinetic description to the hydrodynamic one, which are related and somehow inspired the present work. The gyrokinetic regime, studied by \cite{Golse1999vlasov,Miot2016,saint-raymond_control_2002} provides the convergence of solutions of the two-dimensional Vlasov–Poisson system, in the presence of  very strong magnetic field, to solutions of the incompressible Euler equation in vorticity form. 
The study of the quasineutral regime, starting from seminal works of \cite{brenier_convergence_2000, grenier_oscillations_1996} and also \cite{Golse2003vlasov}, provides hydrodynamic equation in the limit of zero Debye length has now many contributions as \cite{DanielHan2011, han2016quasineutral, griffin-pickering_singular_2020a, porat_derivation_2024}, to name only a few. 
More similar to our result, starting from pressureless Euler systems, the derivation of a number of models has been rigorously established in various dispersive regimes, such as the Zakharov-Kuznetsov equation \cite{lannes_cauchy_2013} and of the KP-II equation \cite{pu_dispersive_2013}. The global dynamics of the ionic Euler-Boltzmann system (with also a Poisson term in the Maxwell-Boltzmann equation) has been addressed in a few works: we refer to \cite{Bae_2024} and \cite{Bae2025} for a discussion on the global well-posedness and blow-up for a 1D model and to the recent work \cite{ertzbischoff2026global}, where the stability of certain profiles, and dispersion phenomena were analyzed. 
Concerning singular limits of hyperbolic systems, we mention \cite{KlainermanMajda1981,klainerman1982compressible} about a general theory for singular limits of quasilinear symmetric hyperbolic systems, proving that solutions converge to a reduced constrained system obtained by projecting the dynamics onto the kernel of the singular operator. As  application, they rigorously justify the low Mach number limit of the compressible Euler equations, showing that, smooth compressible flows converge to solutions of the incompressible Euler equations as the Mach number tends to zero. In \cite{Schochet1994fast}, the author establishes a general theory for fast singular limits of symmetric hyperbolic systems with a large parameter, proving convergence after filtering the fast oscillations  generated by the singular operator. In \cite{Lions1998incompressible}, the authors prove that, under suitable assumptions, global weak solutions of the compressible Navier–Stokes equations converge to a weak solution of the incompressible Navier–Stokes equations as the Mach number tends to zero. In \cite{Feireisl2014multiple}, the authors established the singular limit of the compressible rotating Navier–Stokes equations with general initial data, proving convergence to the two-dimensional incompressible Euler equations by using  relative entropy techniques and dispersive estimates to handle fast acoustic and rotational waves. Finally we mention also the works \cite{bin2012} and \cite{necasova_singular_2020}, closely connected to our, where the author study limit of isothermal Euler equations in the small Mach/Rossby limit for ill-prepared data using completely different techniques. 
\subsubsection*{Structure of the paper}
The rest of the manuscript is organized as follows: in \autoref{Sec:mainresults}, we present  the main results of this work. 
\autoref{Sec:smmoth-solution} is devoted  to the local well posedness of some hyperbolic systems, including \eqref{eul-bol_intro}. In  \autoref{Sec:unif-estimates},  we show the existence  of smooth solution to the singular system \eqref{eul-bol_intro} with a lifespan time uniform in $\varepsilon$. \autoref{Sec:derivtionHM} concerns the identification of the space where the weak limits of \eqref{eul-bol_intro} live, then we introduce the quantity that permit to get strong convergence of an appropriate variable, which allows to get \eqref{HM-eqn-Intro} as a limit of \eqref{eul-bol_intro}. Finally, \autoref{sec:cascade} discuss a  property of
inverse energy cascade and the drift waves for \eqref{HM-eqn-Intro}. 
\section{Main results}\label{Sec:mainresults}
     Let  us  introduce the following subset of $H^s(\mathbb{R}^2), s>2$ 
     \begin{align*}
      \mathbb{B}&=\left\{ V=(\mu,v)\in  (H^s(\mathbb{R}^2))^3;   \text{ and } \Vert V\Vert_{H^s}< +\infty\right\}.
    \end{align*}

Throughout Sections 2--4, except for \autoref{rem:scaling}, the background is a time-independent family $n_0^\varepsilon(x)$, $0<\varepsilon\leq1$, satisfying
\begin{assumption}\label{ass:general-background}
For $s>2$, $n_0^\varepsilon(x)>0$, $\log n_0^\varepsilon \in {H^{s+1}(\mathbb R^2)}$ $ \forall\varepsilon\in (0,1)$ and there exists $g\in H^s_{\operatorname{loc}}(\R^2)$, $g=\nabla G$ such that  
\begin{align*}
  &g_\varepsilon:=\frac1\varepsilon\nabla\log n_0^\varepsilon
 \longrightarrow g\quad\hbox{strongly in }(H_{\operatorname{loc}}^s(\mathbb R^2))^2. 
\\
&\sup_{0<\varepsilon\leq\varepsilon_0}\left(\|g_\varepsilon\|_{L^\infty}+ \|(-\Delta)^{s/2}g_\varepsilon\|_{L^2}\right)< K.
\end{align*}
\end{assumption}
The local well-posedness only needs positivity and $\log n_0^\varepsilon\in H^{s+1}$ for a fixed parameter $0<\varepsilon\leq\varepsilon_0$; the uniform bounds with respect to $\varepsilon$ and the convergence assumptions are used in the subsequent two theorems.
\par

\begin{remark}
    A family satisfying \autoref{ass:general-background} is given by 
    $n_0^\varepsilon(x):= e^{\varepsilon G^\varepsilon(x)}$
    with $G^\varepsilon \in H^{s+1}(\R^2)$, $g_\varepsilon= \nabla G_\varepsilon$ satisfying the above uniform estimates and $\nabla G_\varepsilon \rightarrow \nabla G$ in $(H^s_{\operatorname{loc}}(\R^s))^2$. In particular, we can allow $G^\varepsilon$ to converge to a linear function $v\cdot x$, $v,x\in \R^2$, so that $g$ represents a constant background density gradient $v$. 
\end{remark}

    The first result we prove concern local well-posedness and equivalence of the different reformulations of \eqref{eul-bol_intro}. 

        \begin{theorem}\label{TH1-Main-BIS}  Let $\varepsilon> 0$, $s>2$ and   $n_0^\varepsilon>0$ with $\log n_0^\varepsilon \in H^{s+1}(\R^2)$. Assume that $\big(\phi_0, u_{0}\big)\in \mathbb{B}$:
    \begin{enumerate}
          \item[i)]  There exists  $T_\varepsilon>0$ and a  unique solution $(\phi_\varepsilon,u_\varepsilon)$ to  \eqref{Euler-rho-phi-Main}
         defined on $[0,T_\varepsilon)$ such that
         \begin{align*}
           (\phi_\varepsilon,u_\varepsilon) \in C(0,T_\varepsilon;\mathbb{B}) \cap C^1(0,T_\varepsilon; (H^{s-1}(\mathbb{R}^2))^3),  
         \end{align*}
                  and $T_\varepsilon$ depends only on the initial data, $n_0$ and $\varepsilon$data, $n_0^\varepsilon$ and $\varepsilon$.
            \item[ii)]     Assume that $\rho_0=n_0^\varepsilon e^{\varepsilon\phi_0}>0$ (thus $ \log \rho_0 \in H^s(\R^2) $ and $u_0\in (H^s(\R^2))^2$). Then, there exists a  unique solution $(\rho_\varepsilon,u_\varepsilon,\phi_\varepsilon)$ to  \eqref{eul-bol_intro}  defined on $[0,T_\varepsilon)$,  where $T_\varepsilon$ is the same as above. In particular, we have
         \begin{align*}
u_\varepsilon \in C(0,T_\varepsilon;(H^{s}(\mathbb{R}^2))^2) \cap C^1(0,T_\varepsilon; (H^{s-1}(\mathbb{R}^2))^2)  
         \end{align*}
         and  $\rho_\varepsilon=n_0^\varepsilon e^{\varepsilon\phi_\varepsilon}>0, \  \log \rho_\varepsilon \in C(0,T_\varepsilon;H^s(\R^2))\cap  C^1(0,T_\varepsilon;H^{s-1}(\R^2)).$ Furthermore, 
       we have 
       \begin{align*}
           \rho_\varepsilon \in C(0,T_\varepsilon;C^1_b(\R^2))\cap  C^1(0,T_\varepsilon;C_b(\R^2)).
       \end{align*}
         Finally,  the systems
         \eqref{eul-bol_intro} and \eqref{Euler-rho-phi-Main} are equivalent on  $0<t< T_\varepsilon$.
    \end{enumerate}
    \end{theorem}  

\begin{proof}
    See \autoref{Sec:smmoth-solution}.
\end{proof}
\begin{remark}\label{rem:scaling}
    Consider the following systems for a generic choice of $\tilde n_0>0$.
   \begin{align*}
   \begin{minipage}{0.53\textwidth} 
    \begin{equation}\label{eul-bol_1}
    \begin{cases}
        \partial_t \rho +\Div(\rho u)=0\\
        \partial_t(\rho u) + \Div(\rho u\otimes u)=\rho\left(-\nabla \phi+ u^\perp\right) \\
        \rho(x)=  n_0 e^{\phi}. \tag{EB}
    \end{cases}\end{equation} 
    \end{minipage} 
    \begin{minipage}{0.5\textwidth}
    \begin{equation}\label{PV}\begin{cases}
	&\partial_{t}\phi + \operatorname{div}    u +u \cdot\nabla \left( \phi  + \log n_0 \right)  =0, \\[0.0cm]
&\partial_{t}  u  + u \cdot \nabla u 
=\left(-\nabla\phi
  +  u ^{\perp}\right). \tag{PV}
	\end{cases} \end{equation}
\end{minipage} 
 \end{align*}
where $ \rho(0, x)= \rho_0(x)>0,$ $u(0, x)=u_0(x)$ and $\phi(0,x)= \log(\rho_0(x)/\tilde n_0(x))$. Solutions $(\rho, u)$ and $(\phi, u)$ of the above systems coincide with solutions of \eqref{eul-bol_intro} and \eqref{Euler-rho-phi-Main} for $\varepsilon=1$, respectively. If $s>2$ $\rho_0>0$ and $(\log \rho_0, u_0) \in (H^s(\R^2))^3$ (resp. $(\phi_0, u_0)\in \mathbb{B}$), \autoref{TH1-Main-BIS} ensures that these systems are equivalent.
Now, let  
$$(\phi,u) \in C(0,T_1;\mathbb{B}) \cap C^1(0,T_1; (H^{s-1}(\R^2))^3) $$
be the unique solution of \eqref{PV} associated with the data $(\phi_{0}(x), u_{0}(x))$ and $n_0^ \varepsilon $ given by \autoref{TH1-Main-BIS}.
Introduce the new variables:
\begin{align*}
\phi_{\varepsilon}\left(  t,x\right)    & =\frac{1}{\varepsilon}\phi\left(  \frac{t}{\varepsilon},x\right) ; \quad 
u_{\varepsilon}\left(  t,x\right)     =\frac{1}{\varepsilon}u\left(
\frac{t}{\varepsilon},x\right), \quad   0<t<T_1.
\end{align*}
Then, $(\phi_\varepsilon,u_\varepsilon)  \in C(0,\varepsilon T_1;\mathbb{B}) \cap C^1(0,\varepsilon T_1; (H^{s-1}(\R^2))^3)  
    $ and it is the unique solution to \eqref{Euler-rho-phi-Main} with initial conditions $(\phi^\varepsilon_0, u^\varepsilon_0)= (\varepsilon^{-1}\phi_0, \varepsilon^{-1}u_0)$.
\end{remark}
The next result shows that, if the initial data remain bounded in $\mathbb{B}$, the lifetime of solutions to \eqref{Euler-rho-phi-Main} given by \autoref{TH1-Main-BIS} is actually independent of $\varepsilon$. 
    \begin{theorem}\label{Prop-Hs-estimate}
 Let $s>2$,  let \autoref{ass:general-background} be true and $\big(\phi_0, u_{0}\big)\in \mathbb{B}.$   Let $(\phi_\varepsilon,u_\varepsilon)$ be the unique potential--velocity solution to \eqref{Euler-rho-phi-Main} from \autoref{TH1-Main-BIS}, then
   there exists $T_*>0$, independent of $0<\varepsilon\leq\varepsilon_0$, such that for every $0<T<T_*$ there is $C_T>0$, independent of $\varepsilon$, with
    \begin{align*}
    \sup_{0<t<T}   \| (\phi_\varepsilon(t),u_\varepsilon(t)) \|_{H^s}^2 \leq C_T.
    \end{align*}

Moreover, $T_*$ depends on $y_0:=\|(\phi_0,u_0)\|_{H^s}$ and on $K$, namely
\[
 T_* = \frac{2}{\mathbf M K}\log\left(1+\frac{K}{y_0}\right),
\]
where $\mathbf M>0$ depends only on $s$, the Sobolev embeddings and the Kato--Ponce inequalities. If $K=0$ and $y_0>0$, this formula is understood as $T_*=2/(\mathbf M y_0)$; if $y_0=0$, the solution is identically zero and $T_*=\infty$.
\par
     \end{theorem}
    \begin{proof}
        See \autoref{Sec:unif-estimates}.
    \end{proof}
 
     \begin{remark}
    It is worth mentioning that  the same conclusion of  \autoref{TH1-Main-BIS} and  \autoref{Prop-Hs-estimate} holds in 3D if $s>\dfrac{5}{2}$ and $u^\perp$ is replaced by $u\times e_3$. However, since we are interested in the derivation of Hasegawa-Mima for the 2D system, we formulate all the statement in 2D.
\end{remark} 
Now, we present 
the  main result of this work.
\begin{theorem}[Derivation of Hasegawa-Mima equation]\label{Theorem-main}
Let $s>2$. Let \autoref{ass:general-background} be true and   $\big(\phi_0, u_{0}\big)\in \mathbb{B}$.
    Let $\varepsilon>0$ and $(\phi_\varepsilon,u_\varepsilon)$  
be a  family of solution to \eqref{Euler-rho-phi-Main} associated with the initial data $(\phi_0,u_0).$ Fix $0<T^*<T_*$, with $T_*$ from \autoref{Prop-Hs-estimate}. Then, as $\varepsilon\rightarrow 0$ 
\begin{align*}
&\phi_\varepsilon\buildrel\ast\over\rightharpoonup \phi && \text{in }
L^\infty(0,T^*;H^s(\mathbb R^2))
\\
&u_\varepsilon \buildrel\ast\over\rightharpoonup -\nabla^\perp \phi 
 &&\text{in }L^\infty(0,T^*;(H^s(\mathbb R^2))^2) 
\end{align*}
where $\phi\in  L^\infty(0,T^*;H^{s+1}(\mathbb R^2)) \cap C_w(0,T^*;H^{s+1}(\mathbb{R}^2))$   is the unique  solution to the following Hasegawa-Mima equation:
\begin{align}\label{HM-main-equ}
\tag{HME$_0$}
	\begin{cases}
    & \partial_t(\phi-\Delta\phi)+\nabla^\perp \phi \cdot\nabla( \Delta\phi)-
\nabla^\perp \phi \cdot g=0,\\[0.15cm]
&(\phi-\Delta\phi)\big|_{t=0} = \phi_0 + \nabla^\perp\cdot u_0 \in H^{s-1}(\mathbb{R}^2).
\end{cases}
\end{align}
Moreover, letting $ q_\varepsilon:=\phi_\varepsilon+ \nabla^\perp\cdot u_{\varepsilon}$, we have the following strong convergence, as $\varepsilon\downarrow 0$
\begin{align*}
     q_\varepsilon
\longrightarrow
q:=\phi-\Delta\phi
 \text{  in }
C([0,T^*];H^{s-1-\delta}_{\mathrm{loc}}(\mathbb R^2))  \text{ for any } \delta>0.
\end{align*}

\end{theorem}

The proof will be carried out in \autoref{Sec:derivtionHM}.
\begin{remark}
    By \autoref{ass:general-background}, we have $g= \nabla G$. Thus we get back the equation written as in \eqref{HM-eqn-Intro} just by defining $n_0(x) = e^{G(x)}$.
\end{remark}
\begin{remark}\label{rem:density}
    Note that if one assumes additional to \autoref{ass:general-background} that  $\log n_0^\varepsilon(0)\rightarrow \mathbf{C}$  as $\varepsilon\rightarrow 0$ then $\rho_\varepsilon$ in \eqref{eul-bol_intro} converges to  $e^{\mathbf{C}}$ pointwise, as $\varepsilon\downarrow 0$. Indeed,   recall that 
    \begin{align*}
            \log \rho_\varepsilon(x,t)=\log n_0^\varepsilon(x)+\varepsilon\phi_\varepsilon(t,x), \quad t\in [0,T^*).
    \end{align*}
    The second term vanishes at the limit since $(\phi_\varepsilon)_\varepsilon$ is bounded  thanks to \autoref{Prop-Hs-estimate}. On the other hand, we have
    $$\log n_0^\varepsilon(x) -  \log n_0^\varepsilon(0)= \int_0^1 \frac{d}{ds}\log n_0^\varepsilon(sx)ds. $$
Therefore $|\log n_0^\varepsilon(x) -  \log n_0^\varepsilon(0)| \le|x|\|\nabla\log n_0^\varepsilon\|_{L^\infty} \le \varepsilon K|x|$, thus it converges to zero. Consequently, we get
\begin{align*}
             \lim_{\varepsilon \to 0}   \log \rho_\varepsilon (t,x)=
             \lim_{\varepsilon \to 0}  \log n_0^\varepsilon(x)=
             \lim_{\varepsilon \to 0}  \log n_0^\varepsilon(0)=\mathbf{C}.
\end{align*}
Finally, we remark that \autoref{ass:general-background} forbids $\rho_\varepsilon$ from having finite mass, indeed, if $n_0^\varepsilon$ was integrable, $\log n_0^\varepsilon$ would not be bounded, as requested. 
\end{remark}
\begin{remark} 
It follows from the above theorem that $\partial_t(\phi-\Delta\phi)\in L^\infty(0,T^*;H^{s-2}(\mathbb{R}^2)).$ Moreover, due to the elliptic regularity, the initial condition for $\phi$ satisfies
$$\phi|_{t=0}= (I-\Delta)^{-1}(\phi_0+ \nabla^\perp\cdot u_0) \in H^{s+1}(\mathbb{R}^2).$$ 
This indicate that there is no gain of regularity for the solution of \eqref{HM-main-equ} itself. However, $\phi$ gains one order of regularity compared to $\phi_\varepsilon$ in the limit. 
\end{remark}
\begin{remark}
By exploiting the \(L^2\)-energy estimate for the generalized vorticity $\phi-\Delta\phi$, one can establish the existence of a global weak solution in the same spirit as the analysis of the Hasegawa--Mima equation with singular density in bounded domains \cite{Flandoli-Tahraoui}.


\end{remark}

\section{The scaled potential-velocity system  and  uniform estimates} \label{Sec:smmoth-solution}
\begin{proof}[Proof of \autoref{TH1-Main-BIS}]
Fix $\varepsilon>0$. We prove the theorem directly for this parameter and the background $n_0^\varepsilon$.
Let us start with  system \eqref{Euler-rho-phi-Main}.
   Let $s>2$ and assume that $\log n_0^\varepsilon \in H^{s+1}(\R^2)$ such that   $\big(\phi_0, u_{0}\big)\in \mathbb{B}.$         Set $\Lambda=(\phi_\varepsilon,u_\varepsilon)\in \mathbb{R}^3$ and notice that \eqref{Euler-rho-phi-Main} can be written as follows
    \begin{align}\label{system-2}
        \partial_t \Lambda+\sum_{j=1}^2C_j(\Lambda)\partial_j \Lambda=F_\varepsilon(x,\Lambda), \quad \Lambda\big|_{t=0} =\Lambda_0,
    \end{align}
    where
    \begin{align*}
        C_j(\Lambda)=\begin{pmatrix}
 u_j &  \varepsilon^{-1}e_j^T\\
\varepsilon^{-1}e_j & u_jI_{2\times 2}
\end{pmatrix}, \quad e_1=\begin{pmatrix}
  1\\0
\end{pmatrix}  e_2=\begin{pmatrix}
  0\\1
\end{pmatrix}, F_\varepsilon(x,\Lambda)=\begin{pmatrix}
-u_\varepsilon\cdot g_\varepsilon(x)\\\varepsilon^{-1}u_\varepsilon^\perp
\end{pmatrix}, \quad j=1,2.
    \end{align*}
    Notice that $C_j$ in the system \eqref{system-2} are symmetric.
     Since  $C_j$ in the system \eqref{system-2} is  symmetric, \cite[Theorem 2]{kato1975cauchy}, provides a unique solution $(\phi,u)$ to  \eqref{system-2}
         defined on $[0,T_\varepsilon)$, where $T_\varepsilon>0$, such that
         \begin{align*}
           (\phi_\varepsilon,u_\varepsilon) \in C(0,T_\varepsilon;\mathbb{B}) \cap C^1(0,T_\varepsilon; (H^{s-1}(\mathbb{R}^2))^3)  
         \end{align*}
         and $T_\varepsilon$ depends on the initial data, $n_0^\varepsilon$ and $\varepsilon$.\\
         
Now, we introduce the new variable $\rho$ defined  as $ \rho_\varepsilon(t)=n_0^\varepsilon\exp(\varepsilon\phi_\varepsilon(t))$ for any $t\in [0,T_\varepsilon)$.Since $n_0^\varepsilon>0$ and $\phi  \in H^{s}(\R^2)\hookrightarrow L^\infty(\R^2)$ $\forall t <T_\varepsilon$, one has  $\rho(t) >0$ for any $t\in [0,T_\varepsilon)$.
  On the other hand, since 
    \begin{align*}
\log \rho_\varepsilon=\log n_0^\varepsilon+\varepsilon\phi_\varepsilon \qquad \text{ for any } t\in [0,T_\varepsilon).
    \end{align*}
we have
         \begin{align}
            \log \rho_\varepsilon \in C(0,T_\varepsilon;H^s(\R^2))\cap  C^1(0,T_\varepsilon;H^{s-1}(\R^2)).
         \end{align}
Using the Sobolev embedding $H^s(\R^2) \hookrightarrow C^1_b(\R^2)$ if $s>2$, where  $C^1_b(\R^2)$ is the space of functions whose first-order partial derivatives exist and are continuous and bounded, we get
       \begin{align*}
           \rho_\varepsilon \in C(0,T_\varepsilon;C^1_b(\R^2))\cap  C^1(0,T_\varepsilon;C_b(\R^2)).
       \end{align*}
Finally, by using  $\rho>0$, $\log \rho_\varepsilon=\log n_0^\varepsilon+\varepsilon\phi_\varepsilon$ and \eqref{Euler-rho-phi-Main}, we get that $(\rho_\varepsilon,u_\varepsilon,\phi_\varepsilon)$ is the unique classical solution to  \eqref{eul-bol_intro} (for the fixed $\varepsilon>0$), which completes the proof.
\end{proof}
It is worth mentioning that it is not clear \textit{a priori} whether  $\rho_\varepsilon >0$ if $\rho_0>0$ with general velocity field $u_\varepsilon$, which is important to make sense of \eqref{eul-bol_intro}. However, since  $ u_\varepsilon\in C(0,T^*;H^s(\R^2)), s>2$ and satisfies \eqref{Euler-rho-phi-Main} then a standard computation ensures that 
    $\rho_\varepsilon(t,x)\ge c_{T^*}>0,$
for almost every $(t,x)\in(0,T^*)\times\mathbb{R}^2$, where
\[
c_{T^*}
=
c_0\exp\left(
-\int_0^{T^*}
\|\operatorname{div}u(\tau)\|_{\infty}\,d\tau
\right)
\]
    under the assumption $ \rho_\varepsilon(0,x)=\rho_0(x)\ge c_0>0
\text{ for a.e. }x\in\mathbb{R}^2.$

\subsection{ The scaled  system and  uniform estimates with respect to $\varepsilon$}\label{Sec:unif-estimates}
Thanks to the equivalence provided by \autoref{TH1-Main-BIS}, from now on, we will only deal with \eqref{Euler-rho-phi-Main} to derive the Hasegawa-Mima equations. Recall that the system reads:
        \begin{align}
        \tag{$\varepsilon$-PV}
	\begin{cases}
	&\partial_{t}\phi_\varepsilon + \frac{1}{\varepsilon}\operatorname{div}    u_\varepsilon +u_\varepsilon \cdot \left(\nabla \phi_\varepsilon  + g_\varepsilon \right)  =0, \quad 0<t<T, \quad x\in \mathbb{R}^2,\\[0.0cm]
&\partial_{t}  u_\varepsilon  + u_\varepsilon \cdot \nabla u _\varepsilon   
=\frac{1}{\varepsilon}\left(-\nabla\phi_\varepsilon 
  +  u_\varepsilon ^{\perp}\right).
	\end{cases}
	\end{align} 
    as announced in \autoref{Prop-Hs-estimate}, we now show that the existence time $T_\varepsilon$ for \eqref{Euler-rho-phi-Main} can be chosen uniformly with respect to $0<\varepsilon\leq\varepsilon_0$. The same will then apply also to\eqref{eul-bol_intro} in view of \autoref{TH1-Main-BIS}. \\

 We start by showing the $L^2$-estimate, which indicates that the singular terms does not have contribution in the energy of the system \eqref{Euler-rho-phi-Main}. This observation allows to have a uniform estimate in $H^s, s>2$ in the \autoref{Prop-Hs-estimate}.\\

  Let us test the first equation in \eqref{Euler-rho-phi-Main} by $\phi_\varepsilon$, take the scalar product of the second one by $u_\varepsilon$ and then integrate over $\mathbb R^2$, we get
\begin{align*}
\int_{\mathbb R^2}
\phi_\varepsilon\,\partial_t\phi_\varepsilon\,dx
&
+
\int_{\mathbb R^2}
\phi_\varepsilon(u_\varepsilon\cdot\nabla\phi_\varepsilon)\,dx
\nonumber
\\
&+\frac1\varepsilon
\int_{\mathbb R^2}
\phi_\varepsilon\,\operatorname{div}u_\varepsilon\,dx
+\int_{\mathbb R^2}
\phi_\varepsilon\,u_\varepsilon\cdot g_\varepsilon\,dx
=0,\\
\int_{\mathbb R^2}
u_\varepsilon\cdot\partial_tu_\varepsilon\,dx
&
+
\int_{\mathbb R^2}
u_\varepsilon\cdot(u_\varepsilon\cdot\nabla u_\varepsilon)\,dx
\nonumber
\\
&
+\frac1\varepsilon
\int_{\mathbb R^2}
u_\varepsilon\cdot\nabla\phi_\varepsilon\,dx
-
\frac{1}{\varepsilon}
\int_{\mathbb R^2}
u_\varepsilon^\perp\cdot u_\varepsilon\,dx
=0.
\end{align*}
Let us  estimate each term separately.
First, note that
\[
\int_{\mathbb R^2}
\phi_\varepsilon\,\partial_t\phi_\varepsilon\,dx
=
\frac12
\frac{d}{dt}
\|\phi_\varepsilon\|_{2}^2, \quad \int_{\mathbb R^2}
u_\varepsilon\cdot\partial_tu_\varepsilon\,dx
=
\frac12
\frac{d}{dt}
\|u_\varepsilon\|_{2}^2.
\]
Now, we make the following key observation
\[
\int_{\mathbb R^2}
\phi_\varepsilon\,\operatorname{div}u_\varepsilon\,dx
=
-
\int_{\mathbb R^2}
\nabla\phi_\varepsilon\cdot u_\varepsilon\,dx.
\]
Consequently, we get
\[
\frac1\varepsilon
\int
\phi_\varepsilon\,\operatorname{div}u_\varepsilon
+
\frac1\varepsilon
\int
u_\varepsilon\cdot\nabla\phi_\varepsilon
=
0.
\]
Next, since
$ u^\perp=(u_2,-u_1),
$ we immediately get
$u_\varepsilon^\perp\cdot u_\varepsilon
=0$ and 
$\int_{\mathbb R^2}
u_\varepsilon^\perp\cdot u_\varepsilon\,dx
=0.$
Therefore, singular terms disappear from the energy estimate and we have
\begin{align*}
    \dfrac{1}{2}\frac{d}{dt}(\|\phi_\varepsilon\|_{2}^2+\|u_\varepsilon\|_{2}^2)+
\int_{\mathbb R^2}
\phi_\varepsilon(u_\varepsilon\cdot\nabla\phi_\varepsilon)\,dx
&+ \int_{\mathbb R^2}
u_\varepsilon\cdot(u_\varepsilon\cdot\nabla u_\varepsilon)\,dx\\&=-\int_{\mathbb R^2}
\phi_\varepsilon\,u_\varepsilon\cdot g_\varepsilon\,dx.
\end{align*}
We now estimate the transport terms. The  integration by parts gives
\begin{align*}
\int_{\mathbb R^2}
\phi_\varepsilon(u_\varepsilon\cdot\nabla\phi_\varepsilon)\,dx
&=
\frac12
\int_{\mathbb R^2}
u_\varepsilon\cdot\nabla(\phi_\varepsilon^2)\,dx=
-\frac12
\int_{\mathbb R^2}
(\operatorname{div}u_\varepsilon)\phi_\varepsilon^2\,dx,\\
\int_{\mathbb R^2}
u_\varepsilon\cdot(u_\varepsilon\cdot\nabla u_\varepsilon)\,dx
&=
\frac12
\int_{\mathbb R^2}
u_\varepsilon\cdot\nabla(|u_\varepsilon|^2)\,dx
=
-\frac12
\int_{\mathbb R^2}
(\operatorname{div}u_\varepsilon)|u_\varepsilon|^2\,dx.
\end{align*}
Hence, we obtain
\begin{align*}
   \left|
\int_{\mathbb R^2}
\phi_\varepsilon(u_\varepsilon\cdot\nabla\phi_\varepsilon) dx
\right| 
\le
\frac12
\|\nabla u_\varepsilon\|_{L^\infty}
\|\phi_\varepsilon\|_{2}^2, \quad
\left|
\int_{\mathbb R^2}
u_\varepsilon\cdot(u_\varepsilon\cdot\nabla u_\varepsilon) dx
\right| 
\le
\frac12
\|\nabla u_\varepsilon\|_{L^\infty}
\|u_\varepsilon\|_{2}^2.
\end{align*}

Since
$\|g_\varepsilon\|_{L^\infty}\leq K$, we deduce
\[
\left|
\int_{\mathbb R^2}
\phi_\varepsilon
u_\varepsilon\cdot
 g_\varepsilon
\right| dx
\le
\| g_\varepsilon\|_{L^\infty}
\|\phi_\varepsilon\|_{2}
\|u_\varepsilon\|_{2} \le K \|\phi_\varepsilon\|_{2}
\|u_\varepsilon\|_{2}.
\]
Young's inequality yields
\[
\left|
\int_{\mathbb R^2}
\phi_\varepsilon
u_\varepsilon\cdot
 g_\varepsilon
\right|
\le
\frac12
K
\left(
\|\phi_\varepsilon\|_{2}^2
+
\|u_\varepsilon\|_{2}^2
\right).
\]
Finally,  after collecting the above estimates, we get
\begin{align}\label{L2-esti-varepsilon}
\frac{d}{dt}
\left(
\|\phi_\varepsilon\|_{2}^2
+
\|u_\varepsilon\|_{2}^2
\right)
\le
&
\|\nabla u_\varepsilon\|_{L^\infty}
\left(
\|\phi_\varepsilon\|_{2}^2
+
\|u_\varepsilon\|_{2}^2
\right)
+K
\left(
\|\phi_\varepsilon\|_{2}^2
+
\|u_\varepsilon\|_{2}^2
\right).
\end{align}

We now provide the proof of the main result of this section
\begin{proof}[Proof of \autoref{Prop-Hs-estimate}]\label{Subsection-Hs-est}
The existence and uniqueness part of the statement are provided by \autoref{TH1-Main-BIS}. We are only left to prove an $H^s$-estimate uniform in $0<\varepsilon\leq\varepsilon_0$.
 We recall the following  Bessel potential operator $\Lambda^s =(I-\Delta)^{s/2}$. By applying $\Lambda^s$ to \eqref{Euler-rho-phi-Main}, we get
\begin{align}\label{Euler-rho-phi-Main-Lambda-s}
	\begin{cases}
&\partial_{t}\Lambda^s\phi_{\varepsilon} +\Lambda^s[u_{\varepsilon} \cdot
\nabla\phi_{\varepsilon}]  + \frac{1}{\varepsilon}
\operatorname{div} \Lambda^su_{\varepsilon} +\Lambda^s[u_{\varepsilon
} g_\varepsilon ]     =0,\\[0.2cm]
&\partial_{t}\Lambda^su_{\varepsilon}  +\Lambda^s[
u_{\varepsilon} \cdot\nabla u_{\varepsilon}]  +\frac{1}{\varepsilon}\nabla \Lambda^s\phi_{\varepsilon}-\frac{1}{\varepsilon}
[\Lambda^su_{\varepsilon}]  ^{\perp}   =0.	\end{cases}
	\end{align}
 Let us test the first equation in \eqref{Euler-rho-phi-Main-Lambda-s} by $\Lambda^s\phi_\varepsilon$, take the scalar product of the second one by $\Lambda^su_\varepsilon$ and then integrate over $\mathbb R^2$, we get
\begin{align*}
\int_{\mathbb R^2}\partial_{t}\Lambda^s\phi_{\varepsilon}  \Lambda^s\phi_{\varepsilon}dx&+ \int_{\mathbb R^2}\Lambda^s[u_{\varepsilon} \cdot
\nabla\phi_{\varepsilon}]  \Lambda^s\phi_{\varepsilon}  dx \\&+ \frac{1}{\varepsilon}
\int_{\mathbb R^2}\operatorname{div} (\Lambda^su_{\varepsilon}) \Lambda^s\phi_{\varepsilon}  dx +\int_{\mathbb R^2}\Lambda^s[u_{\varepsilon
}  \cdot  g_\varepsilon ]   \Lambda^s\phi_{\varepsilon}  dx   =0,\\[0.2cm]
\int_{\mathbb R^2}\partial_{t}\Lambda^su_{\varepsilon}\cdot \Lambda^su_{\varepsilon} dx  &+\int_{\mathbb R^2}\Lambda^s[
u_{\varepsilon} \cdot\nabla u_{\varepsilon}] \cdot \Lambda^su_{\varepsilon} dx  \\&+\frac{1}{\varepsilon}\int_{\mathbb R^2}\nabla \Lambda^s\phi_{\varepsilon}\cdot \Lambda^su_{\varepsilon} dx- \frac{1}{\varepsilon}\int_{\mathbb R^2}
[\Lambda^su_{\varepsilon}]  ^{\perp}  \cdot \Lambda^su_{\varepsilon} dx  =0.
\end{align*}
Note that
\begin{align*}
    \int_{\mathbb R^2}
\partial_t\Lambda^s\phi_\varepsilon\,\Lambda^s\phi_\varepsilon\,dx
=
\frac12
\frac{d}{dt}
\|\Lambda^s\phi_\varepsilon\|_{2}^2, \quad \int_{\mathbb R^2}
\partial_t\Lambda^su_\varepsilon\cdot\Lambda^su_\varepsilon\,dx
=
\frac12
\frac{d}{dt}
\|\Lambda^su_\varepsilon\|_{2}^2.
\end{align*}
As with the $L^2$-estimate, singular terms disappear from the evolution of $H^s$-norm. Indeed, after integration by parts, we obtain
\[
\int_{\mathbb R^2}
\operatorname{div}(\Lambda^su_\varepsilon) \Lambda^s\phi_\varepsilon\,dx
=
-
\int_{\mathbb R^2}
\nabla \Lambda^s\phi_\varepsilon\cdot \Lambda^su_\varepsilon\,dx.
\]
Consequently, we get
\[
\frac1\varepsilon
\int
\operatorname{div}(\Lambda^su_\varepsilon) \Lambda^s\phi_\varepsilon\ dx
+
\frac1\varepsilon
\int
\nabla \Lambda^s\phi_\varepsilon \cdot\Lambda^su_\varepsilon dx
=
0.
\]
Next, since
$ u^\perp=(u_2,-u_1),
$ we immediately get
$\int_{\mathbb R^2}
[\Lambda^su_{\varepsilon}]  ^{\perp}  \cdot \Lambda^su_{\varepsilon} dx
=0.$
Therefore,   we have
\begin{align*}
    \dfrac{1}{2}\frac{d}{dt}(\|\Lambda^s\phi_\varepsilon\|_{2}^2&+\|\Lambda^su_\varepsilon\|_{2}^2)+\int_{\mathbb R^2}\Lambda^s[u_{\varepsilon} \cdot
\nabla\phi_{\varepsilon}]  \Lambda^s\phi_{\varepsilon}  dx
+ \int_{\mathbb R^2}\Lambda^s[
u_{\varepsilon} \cdot\nabla u_{\varepsilon}] \cdot \Lambda^su_{\varepsilon} dx\\&=- \int_{\mathbb R^2}\Lambda^s[u_{\varepsilon
}  \cdot g_\varepsilon ]\Lambda^s\phi_\varepsilon dx.
\end{align*}

Let us estimate the transport terms. Note that 
\begin{align*}
    &\int_{\mathbb R^2}\Lambda^s[u_{\varepsilon} \cdot
\nabla\phi_{\varepsilon}]  \Lambda^s\phi_{\varepsilon}  dx\\&=  \int_{\mathbb R^2}u_{\varepsilon} \cdot
\nabla \Lambda^s\phi_{\varepsilon} \Lambda^s\phi_{\varepsilon}  dx+ \int_{\mathbb R^2} \Bigg[ \Lambda^s[u_{\varepsilon} \cdot
\nabla\phi_{\varepsilon}]  \Lambda^s\phi_{\varepsilon}-u_{\varepsilon} \cdot
\nabla \Lambda^s\phi_{\varepsilon} \Lambda^s\phi_{\varepsilon}\Bigg] dx,\\
&\int_{\mathbb R^2}\Lambda^s[
u_{\varepsilon} \cdot\nabla u_{\varepsilon}] \cdot \Lambda^su_{\varepsilon} dx\\&= \int_{\mathbb R^2}
u_{\varepsilon} \cdot \nabla \Lambda^su_{\varepsilon} \cdot \Lambda^su_{\varepsilon} dx+  \int_{\mathbb R^2} \Bigg[\Lambda^s[
u_{\varepsilon} \cdot\nabla u_{\varepsilon}] \cdot \Lambda^su_{\varepsilon}-u_{\varepsilon} \cdot \nabla\Lambda^s u_{\varepsilon} \cdot \Lambda^su_{\varepsilon} \Bigg] dx.
\end{align*}
As with the $L^2$ estimate, we can bound the first terms of the right-hand side of the last equalities. Integration by parts gives
\begin{align*}
    \int_{\mathbb R^2}
(u_\varepsilon\cdot\nabla \Lambda^s\phi_\varepsilon)\Lambda^s\phi_\varepsilon\,dx
&=
\frac12
\int_{\mathbb R^2}
u_\varepsilon\cdot\nabla( (\Lambda^s\phi_\varepsilon)^2)\,dx=
-\frac12
\int_{\mathbb R^2}
(\operatorname{div}u_\varepsilon) \Lambda^s\phi_\varepsilon^2\,dx,\\
\int_{\mathbb R^2}
(u_\varepsilon\cdot\nabla \Lambda^su_\varepsilon) \cdot  \Lambda^su_\varepsilon\,dx
&=
\frac12
\int_{\mathbb R^2}
u_\varepsilon\cdot\nabla(|\Lambda^su_\varepsilon|^2)\,dx
=
-\frac12
\int_{\mathbb R^2}
(\operatorname{div}u_\varepsilon)|\Lambda^su_\varepsilon|^2\,dx.
\end{align*}
Hence, we obtain
\begin{align*}
 &  \left|
\ \int_{\mathbb R^2}
(u_\varepsilon\cdot\nabla \Lambda^s\phi_\varepsilon)\Lambda^s\phi_\varepsilon\,dx
\right| 
\le
\frac12
\|\nabla u_\varepsilon\|_{L^\infty}
\|\Lambda^s\phi_\varepsilon\|_{2}^2,\\ 
&\left|
\int_{\mathbb R^2}
(u_\varepsilon\cdot\nabla \Lambda^su_\varepsilon) \cdot  \Lambda^su_\varepsilon\,dx
\right|
\le
\frac12
\|\nabla u_\varepsilon\|_{L^\infty}
\|\Lambda^su_\varepsilon\|_{2}^2.
\end{align*}
Concerning the remaining parts we use Cauchy Schwarz  inequality to get
\begin{align*}
  \left|  \int_{\mathbb R^2} \Bigg[ \Lambda^s[u_{\varepsilon} \cdot
\nabla\phi_{\varepsilon}]  \Lambda^s\phi_{\varepsilon}-u_{\varepsilon} \cdot
\nabla \Lambda^s\phi_{\varepsilon} \Lambda^s\phi_{\varepsilon}\Bigg] dx\right| \leq \Vert \Lambda^s[u_{\varepsilon} \cdot
\nabla\phi_{\varepsilon}]-u_{\varepsilon} \cdot
\nabla \Lambda^s\phi_{\varepsilon}\Vert_2\Vert \Lambda^s\phi_{\varepsilon}\Vert_2.
\end{align*}
Now, we use  the Kato--Ponce commutator inequality (see \cite[Lemma X1]{kato1988commutator}) to get
\begin{align*}
    \Vert \Lambda^s[u_{\varepsilon} \cdot
\nabla\phi_{\varepsilon}]-u_{\varepsilon} \cdot
\nabla \Lambda^s\phi_{\varepsilon}\Vert_2 \leq  C(\Vert \nabla u_{\varepsilon}\Vert_\infty\Vert \Lambda^{s-1}\nabla \phi_{\varepsilon} \Vert_2 +\Vert\Lambda^s u_{\varepsilon} \Vert_2\Vert \nabla\phi_\varepsilon\Vert_\infty)
\end{align*}
where $C>0$ is independent of $\varepsilon$ (here and later, we denote with $C$ a positive constants  that can change from one line to the other).  Since \(s>2\), Sobolev embedding ensures
\begin{align*}
    \Vert \nabla u_{\varepsilon}\Vert_\infty \leq M\Vert  u_{\varepsilon}\Vert_{H^s} \text{ and }     \Vert \nabla \phi_{\varepsilon}\Vert_\infty \leq M \Vert \phi_{\varepsilon}\Vert_{H^s},
\end{align*}
where $M>0$ depends only on the embedding $H^s(\mathbb R^2)\hookrightarrow W^{1,\infty}(\mathbb R^2).$
Consequently, there exists $C>0$ independent of $\varepsilon$ such that
\begin{align*}
     \Vert \Lambda^s[u_{\varepsilon} \cdot
\nabla\phi_{\varepsilon}]-u_{\varepsilon} \cdot
\nabla \Lambda^s\phi_{\varepsilon}\Vert_2 \leq C \Vert \phi_{\varepsilon}\Vert_{H^s}\Vert  u_{\varepsilon}\Vert_{H^s}.
\end{align*}
By repeating the same argument, we obtain
\begin{align*}
  \left|  \int_{\mathbb R^2} \Bigg[\Lambda^s[
u_{\varepsilon} \cdot\nabla u_{\varepsilon}] \cdot \Lambda^su_{\varepsilon}-u_{\varepsilon} \cdot \nabla\Lambda^s u_{\varepsilon} \cdot \Lambda^su_{\varepsilon} \Bigg] dx \right| \leq  \Vert \Lambda^s[
u_{\varepsilon} \cdot\nabla u_{\varepsilon}] -u_{\varepsilon} \cdot \nabla\Lambda^s u_{\varepsilon}\Vert_2\Vert \Lambda^s u_{\varepsilon}\Vert_2,
\end{align*}
and 
\begin{align*}
    \Vert \Lambda^s[
u_{\varepsilon} \cdot\nabla u_{\varepsilon}] -u_{\varepsilon} \cdot \nabla\Lambda^s u_{\varepsilon}\Vert_2 \leq C \Vert  u_{\varepsilon}\Vert_{H^s}^2.
\end{align*}
Finally, there exists $C>0$ independent of $\varepsilon$ such that
\begin{align*}
  \left|  \int_{\mathbb R^2}\Lambda^s[u_{\varepsilon
}  \cdot  g_\varepsilon ]\Lambda^s\phi_\varepsilon dx\right| \leq CK\Vert u_\varepsilon \Vert_{H^s}\Vert \phi_\varepsilon \Vert_{H^s}.
\end{align*}
Indeed,  we have
\begin{align*}
     \left|  \int_{\mathbb R^2}\Lambda^s[u_{\varepsilon
}  \cdot g_\varepsilon ]\Lambda^s\phi_\varepsilon dx\right| \leq \Vert \Lambda^s[u_{\varepsilon
}  \cdot g_\varepsilon ]\Vert_2\Vert \Lambda^s\phi_\varepsilon\Vert_2.
\end{align*}

Denote by $D^s=(-\Delta)^{s/2}$. The Kato--Ponce inequality \cite[Theorem 1]{Grafakos03062014} gives
\begin{align*}
 \|\Lambda^s[u_\varepsilon\cdot g_\varepsilon]\|_2
 &\leq C\big(\|D^su_\varepsilon\|_2\|g_\varepsilon\|_\infty
 +\|D^sg_\varepsilon\|_2\|u_\varepsilon\|_\infty
 +\|g_\varepsilon\|_\infty\|u_\varepsilon\|_2\big)\\
 &\leq CK\|u_\varepsilon\|_{H^s}.
\end{align*}
Here $\|g_\varepsilon\|_{H^s}\leq K$ and $H^s(\mathbb R^2)\hookrightarrow L^\infty(\mathbb R^2)$.
Gathering the above estimates, we get
\begin{align*}
    \dfrac{1}{2}\frac{d}{dt}(\|\Lambda^s\phi_\varepsilon\|_{2}^2&+\|\Lambda^su_\varepsilon\|_{2}^2)\leq  
\frac12
\|\nabla u_\varepsilon\|_{L^\infty}( \|\Lambda^s\phi_\varepsilon\|_{2}^2+
\|\Lambda^su_\varepsilon\|_{2}^2)+C \Vert  u_{\varepsilon}\Vert_{H^s}[\Vert  u_{\varepsilon}\Vert_{H^s}^2+ \Vert \phi_{\varepsilon}\Vert_{H^s}^2]\\
&+CK\Vert u_\varepsilon \Vert_{H^s}\Vert \phi_\varepsilon \Vert_{H^s}.
\end{align*}
Since $H^s(\mathbb R^2)\hookrightarrow W^{1,\infty}(\mathbb R^2)$ and $\|\Lambda^s\phi_\varepsilon\|_{2}^2+\|\Lambda^su_\varepsilon\|_{2}^2$ is equivalent norm to   $\| (\phi_\varepsilon,u_\varepsilon) \|_{H^s}^2 $, we deduce
\begin{align}\label{inequality-est-local}
    \frac{d}{dt}\| (\phi_\varepsilon,u_\varepsilon) \|_{H^s}^2&\leq  C \Vert  u_{\varepsilon}\Vert_{H^s}\|(\phi_\varepsilon,u_\varepsilon) \|_{H^s}^2+CK\|(\phi_\varepsilon,u_\varepsilon) \|_{H^s}^2\\
    & \notag\leq C\|(\phi_\varepsilon,u_\varepsilon) \|_{H^s}^3+C_2\|(\phi_\varepsilon,u_\varepsilon) \|_{H^s}^2.
\end{align}
If we denote by $f(t)= \Vert  (\phi_\varepsilon(t),u_\varepsilon(t))\Vert_{H^s}^2$, we have
\begin{align*}
    f^\prime (t) \leq Cf^{3/2}(t)+C_2 f(t), \quad f(0)= \Vert  (\phi_0,u_0)\Vert_{H^s}^2.
\end{align*}
As a consequence, the solution stays bounded on $[0,T]$ for every $T<T_*$, where
\[
T_*=
\displaystyle
\frac{2}{C_2}\log\!\left(
1+\frac{C_2}{C\,\|(\phi_0,u_0)\|_{H^s}}
\right).
\]

The conventions for $K=0$ and zero initial data are as in \autoref{Prop-Hs-estimate}. For every fixed $T<T_*$, this differential inequality gives a bound $C_T$ independent of $0<\varepsilon\leq\varepsilon_0$. The continuation criterion for the symmetric hyperbolic system then extends each solution through $T$, proving the asserted uniform lifetime.
\end{proof}

\begin{remark}
We perform the above computations directly using integration by parts and commutation between $\Lambda^s$ and derivatives. This is true in the class of  smooth function. In the previous computations, we focused on the evolution of the $H^s$-norm  of the unique solution to \eqref{Euler-rho-phi-Main}. We can therefore consider the unique solution of the same system with \textit{e.g.} $H^{s+2}$-data , where all computations are justified and the final estimate about the $H^s$-norm can be obtained as  a limit by approximating the data. We have omitted this detail here since it is standard.
\end{remark}
\section{Derivation of the Hasegawa-Mima equation}\label{Sec:derivtionHM}
Before deriving the Hasegawa-Mima equation, let us begin with the following remarks which give an indication of the structure of the limit of \eqref{Euler-rho-phi-Main}.\\

From \autoref{Prop-Hs-estimate}, we have a uniform estimate with respect to $\varepsilon$. However, it is clear from the system 
\eqref{Euler-rho-phi-Main} that $\partial_t (\phi_\varepsilon,u_\varepsilon)$ is of order $ {\varepsilon}^{-1}$ for every $0<\varepsilon\leq\varepsilon_0$ in appropriate topology and getting, \textit{a priori}, a uniform estimate for the $(\partial_t (\phi _\varepsilon,u_\varepsilon))_\varepsilon$ does not seem feasible.  On the other hand, the presence of the singular terms of order $\frac{1}{\varepsilon}$ imposes a
constraint on any possible weak limit of the sequence
$Z_\varepsilon=(\phi_\varepsilon,u_\varepsilon)$. Indeed, the system can be
rewritten in the abstract form
\[
\partial_t Z_\varepsilon+\mathcal{N}_\varepsilon(Z_\varepsilon)
+\frac{1}{\varepsilon}\mathcal{L}(Z_\varepsilon)=0,
\]
where  the singular operator
 $\mathcal{L}$  is given by
\[
\mathcal{L}(\phi,u)=\left(\operatorname{div}u,\nabla\phi-u^\perp\right)
\]
and $\mathcal{N}_\varepsilon$ contains the  remaining terms. Assume that, up to the extraction of a subsequence,
$Z_\varepsilon\rightharpoonup Z=(\phi,u)$ in $L^\infty(0,T^*;(H^s(\R^d))^3$.
Since  the
singular contribution must vanish in the limit in order to prevent the
appearance of unbounded oscillations, one obtains
$
\mathcal{L}(Z_\varepsilon)\longrightarrow 0
$
in the sense of distributions. Passing to the limit yields
$\mathcal{L}(Z)=0,$
which means that every weak limit belongs to the kernel of the singular
operator:
\[
\ker\mathcal{L}
=
\left\{(\phi,u):\operatorname{div}u=0,\quad
\nabla\phi=u^\perp\right\}.
\]
In other words, the singular part of the dynamics acts as a constraint operator: the fast oscillations generated by $\mathcal{L}$ vanish in the limit, and the
limiting state evolves inside the manifold of equilibria of the fast system. Based on  this observation,  let us  show how to obtain the Hasegawa-Mima equation from \eqref{Euler-rho-phi-Main}.
For that, we proceed in many steps.

\subsection{Step 1: the weak limits}
     \autoref{Prop-Hs-estimate}, gives a family of solutions to \eqref{Euler-rho-phi-Main}, that satisfy
    \begin{align*}
        \sup_{0<\varepsilon\leq\varepsilon_0}
\left(
\|\phi_\varepsilon\|_{L^\infty(0,T^*;H^s)}
+
\|u_\varepsilon\|_{L^\infty(0,T^*;H^s)}
\right)
\leq C,
\qquad s>2.
    \end{align*}
    The above uniform estimate give, by the
Banach--Alaoglu theorem, that there exists
$(\phi,u)\in L^\infty(0,T^*;H^s(\mathbb R^2))\times L^\infty(0,T^*;(H^s(\mathbb R^2))^2)$
and a subsequence, still denoted by $\varepsilon$, such that
\begin{align}\label{weak-limi-cv}
    \phi_\varepsilon\buildrel\ast\over\rightharpoonup \phi,
\text{ and }
u_\varepsilon \buildrel\ast\over\rightharpoonup u
\quad\text{in }
L^\infty(0,T^*;H^s(\mathbb R^2))  \text{ and }L^\infty(0,T^*;(H^s(\mathbb R^2))^2) \text{ respectively.}
\end{align}
   
\subsection{Step 2: identification of the weak limits}
\begin{lemma}\label{lem-limits}
Let $0<\varepsilon\leq\varepsilon_0$ and $(\phi_\varepsilon,u_\varepsilon)$ be a family of  solutions to \eqref{Euler-rho-phi-Main}, given by \autoref{TH1-Main-BIS}, and $(\phi,u)$ given by \eqref{weak-limi-cv}.
Then,  the limit satisfies
$
\operatorname{div}u=0,$ $
u^\perp=\nabla\phi
$
and consequently
\[
u=-\nabla^\perp\phi .
\]
\end{lemma}
\begin{proof}
First, we  prove that the limit velocity is divergence free. From the first
equation of \eqref{Euler-rho-phi-Main}, we have
\[
\operatorname{div}u_\varepsilon
=
-\varepsilon\partial_t\phi_\varepsilon
-\varepsilon u_\varepsilon\cdot\nabla\phi_\varepsilon
-\varepsilon
u_\varepsilon\cdot g_\varepsilon.
\]

Let $\varphi\in C_c^\infty((0,T^*)\times\mathbb R^2)$. Then
\begin{align}
 \label{eq-div-1} \left\langle
\operatorname{div}u_\varepsilon,\varphi
\right\rangle
=
&\varepsilon
\int_0^{T^*}\int_{\mathbb R^2}
\phi_\varepsilon\partial_t\varphi\,dxdt
-\varepsilon
\int_0^{T^*}\int_{\mathbb R^2}
(u_\varepsilon\cdot\nabla\phi_\varepsilon)\varphi\,dxdt\\&-\varepsilon
\int_0^{T^*}\int_{\mathbb R^2}
u_\varepsilon\cdot
 g_\varepsilon
\varphi\,dxdt .   \notag
\end{align}
Since $s>2$,  $u_\varepsilon\cdot\nabla\phi_\varepsilon \in H^{s-1}(\R^2)$ (see  \cite[Lemma 2.1.]{kato1975cauchy}
and therefore
$
\|u_\varepsilon\cdot\nabla\phi_\varepsilon\|_
{L^\infty(0,T^*;H^{s-1})}
\leq C.
$
Moreover, we have
$\| g_\varepsilon\|_{\infty}
\leq
K.$
Consequently, all the terms on the right-hand side in \eqref{eq-div-1}  converge to zero
as $\varepsilon\downarrow 0$. Hence, we get
\[
\operatorname{div}u_\varepsilon\longrightarrow 0
\quad\text{in }
\mathcal D'((0,T^*)\times\mathbb R^2).
\]
On the other hand, passing to the limit using the weak-* convergence \eqref{weak-limi-cv} of $u_\varepsilon$ gives
\[
\operatorname{div}u_\varepsilon \to \operatorname{div}u=0
\quad\text{in }
\mathcal D'((0,T^*)\times\mathbb R^2).
\]

Next, we  identify the relation between $\phi$ and $u$. The second equation of
\eqref{Euler-rho-phi-Main} gives
\[
\nabla\phi_\varepsilon-u_\varepsilon^\perp
=
-\varepsilon
\left(
\partial_tu_\varepsilon
+
u_\varepsilon\cdot\nabla u_\varepsilon
\right).
\]
Let
$\psi\in C_c^\infty((0,T^*)\times\mathbb R^2;\mathbb R^2).$
Then
\begin{align}\label{weak-cv-2-orto}
    \left\langle
\nabla\phi_\varepsilon-u_\varepsilon^\perp,\psi
\right\rangle
=
\varepsilon
\int_0^{T^*}\int_{\mathbb R^2}
u_\varepsilon\cdot\partial_t\psi\,dxdt-\varepsilon
\int_0^{T^*}\int_{\mathbb R^2}
(u_\varepsilon\cdot\nabla u_\varepsilon)
\cdot\psi\,dxdt .
\end{align}

The first term  in \eqref{weak-cv-2-orto} converges to zero due to the uniform bound of
$u_\varepsilon$ in $L^\infty(0,T^*;H^s)$. Furthermore, again by using \cite[Lemma 2.1.]{kato1975cauchy}, we have
\[
\|u_\varepsilon\cdot\nabla u_\varepsilon\|_
{L^\infty(0,T^*;H^{s-1})}
\leq C,
\]
so the second term  in \eqref{weak-cv-2-orto} also converges to zero. Therefore, we deduce
\[
\nabla\phi_\varepsilon-u_\varepsilon^\perp
\longrightarrow 0
\quad\text{in }
\mathcal D'((0,T^*)\times\mathbb R^2).
\]

On the other hand, the weak-* convergence \eqref{weak-limi-cv} implies
\[
\nabla\phi_\varepsilon\rightharpoonup\nabla\phi,
\qquad
u_\varepsilon^\perp\rightharpoonup u^\perp \text{ in  } D'((0,T^*)\times\mathbb R^2).
\]
Passing to the limit yields
\[
 \nabla\phi_\varepsilon-u_\varepsilon^\perp \to  \nabla\phi-u^\perp=0
\quad\text{in }
\mathcal D'((0,T^*)\times\mathbb R^2).
\]
Hence, $u^\perp=\nabla\phi .$
Finally, since
$
(V^\perp)^\perp=-V
$
for vectors in $\mathbb R^2$, we apply the perpendicular operator to the
relation above and we get
$
u=-\nabla^\perp\phi .
$
Therefore every weak limit $(\phi,u)$ satisfies
\begin{align}\label{Ker-limit-oper}
    \operatorname{div}u=0,
\qquad
u=-\nabla^\perp\phi.
\end{align}
\end{proof}

    We conclude this part by the following remark.
\begin{remark}\label{Rmq-regularity-limit} Since $u\in L^\infty(0,T^*;(H^s(\mathbb R^2))^2) $ and $\phi \in L^\infty(0,T^*;H^s(\mathbb R^2))$, the last equality \eqref{Ker-limit-oper} ensures that $\phi \in L^\infty(0,T^*;H^{s+1}(\mathbb R^2)),$ which gives a gain in regularity  for the limit of the system \eqref{Euler-rho-phi-Main}. However, any weak limit lives in infinite-dimensional space $\mathcal{A}$ given by
 \begin{align}
    \mathcal{A}=\{ (\phi, -\nabla^\perp \phi ), \quad \phi \in   L^\infty(0,T^*;H^{s+1}(\mathbb R^2))\}.
\end{align}
\end{remark}
 
Our aim, in the next subsection, is to show that $\phi$ in fact is the unique smooth solution (see \autoref{cor-limit-HM}) of the following Hasegawa-Mima equation:
\begin{align*}
	\begin{cases}
    & \partial_t(\phi-\Delta\phi)+\nabla^\perp \phi \cdot\nabla( \Delta\phi)-
\nabla^\perp \phi \cdot g=0,\\[0.15cm]
&(\phi-\Delta\phi)\big|_{t=0} = \phi_0+\nabla^\perp\cdot u_0 \in H^{s-1}(\mathbb{R}^2).
\end{cases}
\end{align*}
\subsection{Step 3: vorticity form and cancellation of the fast oscillation}\label{subsec:vortcity}
Let us  introduce the scalar vorticity:
$\omega_{\varepsilon}:=-\nabla^\perp\cdot  u_\varepsilon
=\partial_1 u_{\varepsilon,2}
-\partial_2 u_{\varepsilon,1}.$
Taking the curl of the momentum equation in  \eqref{Euler-rho-phi-Main} and 
using the identity $-\nabla^\perp\cdot(\nabla\phi_\varepsilon)=0,$ we get
\[
\partial_t\omega_\varepsilon
-\nabla^\perp\cdot 
(u_\varepsilon\cdot\nabla u_\varepsilon)
+\frac{1}{\varepsilon}
\nabla^\perp\cdot (u_\varepsilon^\perp)=0.
\]
Now, note that 
$
-\nabla^\perp\cdot (u_\varepsilon^\perp)
=
-\nabla^\perp\cdot (u_{\varepsilon,2},-u_{\varepsilon,1})
=
\partial_1(-u_{\varepsilon,1})
-\partial_2(u_{\varepsilon,2})
$
and therefore
$$-\nabla^\perp\cdot(u_\varepsilon^\perp)
=
-\operatorname{div}u_\varepsilon.$$
Since $u_\varepsilon$ is two-dimensional, we have 
$
-\nabla^\perp\cdot
(u_\varepsilon\cdot\nabla u_\varepsilon)
=
u_\varepsilon\cdot\nabla\omega_\varepsilon
+\omega_\varepsilon\operatorname{div}u_\varepsilon .
$
Hence the vorticity equation is given by
\begin{align}\label{vorticity}
  \partial_t\omega_\varepsilon
+u_\varepsilon\cdot\nabla\omega_\varepsilon
+\omega_\varepsilon\operatorname{div}u_\varepsilon
+\frac{1}{\varepsilon}\operatorname{div}u_\varepsilon=0 .  
\end{align}
By subtracting \eqref{vorticity} from the first equation  of \eqref{Euler-rho-phi-Main}, we write
\begin{align}\label{equ-derive-HM}
	\begin{cases}
    &\partial_t q_\varepsilon
+
u_\varepsilon\cdot\nabla q_\varepsilon
-\omega_\varepsilon
\operatorname{div}u_\varepsilon
+
u_\varepsilon\cdot g_\varepsilon
=0,\\[0.15cm]
&q_\varepsilon=
\phi_\varepsilon+ \nabla^\perp\cdot u_{\varepsilon},\\[0.15cm]
&\big(\phi_\varepsilon, u_\varepsilon\big)    \big|_{t=0} = \big(\phi_0, u_{0}\big).
\end{cases}
\end{align}
Note that $q_\varepsilon$ 
is the unique solution to \eqref{equ-derive-HM} and satisfies \[ q_\varepsilon=\phi_\varepsilon +  \nabla^\perp\cdot u_{\varepsilon} \in C(0,T^*;H^{s-1}(\mathbb{R}^2))
\cap C^1(0,T^*;H^{s-2}(\mathbb{R}^2))
 \text{ for every } \varepsilon>0.
\] Moreover, the following uniform estimate holds.
\begin{lemma}\label{Lem-bound-q-appro}  Let $s>2,$ then $(q_\varepsilon)_\varepsilon$ is uniformly bounded in  $C(0,T^*;H^{s-1}(\mathbb{R}^2)) \cap C^1(0,T^*;H^{s-2}(\mathbb{R}^2))$.
\end{lemma}
\begin{proof}
    The boundedness in  $C(0,T^*;H^{s-1}(\mathbb{R}^2))$ is a direct consequence of  \autoref{Prop-Hs-estimate} and the definition of $q_\varepsilon$. On the other hand, from the first equation in \eqref{equ-derive-HM}, we have
    \begin{align}\label{time-deri-est}
        \partial_t q_\varepsilon
=
-u_\varepsilon\cdot\nabla q_\varepsilon
+\omega_\varepsilon
\operatorname{div}u_\varepsilon
-
u_\varepsilon\cdot g_\varepsilon.
    \end{align}
    Thanks to \cite[Lemma 2.1.]{kato1975cauchy}, we have
    \begin{align}\label{est-der-1}
    \Vert u_\varepsilon\cdot\nabla q_\varepsilon   \Vert_{H^{s-2}} \leq      \Vert u_\varepsilon\Vert_{H^{s}}    \Vert \nabla q_\varepsilon   \Vert_{H^{s-2}} \text{ and }      \Vert \omega_\varepsilon
\operatorname{div}u_\varepsilon \Vert_{H^{s-2}} \leq      \Vert \operatorname{div}u_\varepsilon\Vert_{H^{s-1}}    \Vert \omega_\varepsilon   \Vert_{H^{s-1}}.
    \end{align}
An argument already detailed guarantees
    \begin{align}\label{est-der-2}
    \Vert u_\varepsilon\cdot g_\varepsilon\Vert_{H^{s}}  \leq    C K\Vert u_\varepsilon \Vert_{H^s}.
    \end{align}
    Since $H^s(\R^s) \hookrightarrow  H^{s-2}(\R^2),$ \eqref{time-deri-est} combined with \eqref{est-der-1}, \eqref{est-der-2} and \autoref{Prop-Hs-estimate} ensures that $( \partial_t q_\varepsilon)_\varepsilon$ is uniformly bounded in  $C(0,T^*;H^{s-2}(\mathbb{R}^2))$, which completes the proof.
\end{proof}
Before presenting a result on the strong convergence of $(q_\varepsilon)_\varepsilon$ and the derivation of Hasegawa-Mima equation, we have the following convergence results.
\begin{lemma}\label{Lem-log} Under 
\autoref{ass:general-background} with $s>2$,
and that
\[
u_\varepsilon\buildrel\ast\over\rightharpoonup u
\quad\text{in }
L^\infty(0,T^*;H^s(\mathbb R^2)).
\]
It holds
$
u_\varepsilon\cdot g_\varepsilon
\buildrel\ast\over\rightharpoonup
u\cdot g \text{ in } L^\infty(0,T^*;H^{s}(\mathbb R^2))
 \text{ as } \varepsilon\downarrow 0. $
\end{lemma}
\begin{proof}
First,  from \eqref{est-der-2} recall that $(u_\varepsilon\cdot g_\varepsilon)_\varepsilon$ is bounded in  $L^\infty(0,T^*;H^{s}(\mathbb R^2)).$ By  the Banach–Alaoglu
theorem, it converges weakly-* (up to a subsequence) to some limit $\xi\in L^\infty(0,T^*;H^{s}(\mathbb R^2)).$
We show that $\xi=u\cdot g$. Indeed, 
let
$
\varphi\in C_c^\infty((0,T^*)\times\mathbb R^2),
$
we write
\begin{align*}
    &\int_0^{T^*}\int_{\mathbb R^2}
u_\varepsilon\cdot g_\varepsilon
\varphi\,dxdt
-
\int_0^{T^*}\int_{\mathbb R^2}
u\cdot g\varphi\,dxdt
\\
&=
\int_0^{T^*}\int_{\mathbb R^2}
(u_\varepsilon-u)\cdot
g\varphi\,dxdt
+
\int_0^{T^*}\int_{\mathbb R^2}
u_\varepsilon\cdot
\left( g_\varepsilon- g
\right)
\varphi\,dxdt.
\end{align*}
The first term converges to zero since 
$u_\varepsilon\buildrel\ast\over\rightharpoonup u
\quad\text{in }
L^\infty(0,T^*;H^s).
$
For the second term, the compactness of the support of $\varphi$, the convergence
\[
\|g_\varepsilon-g\|_{H^s_{\operatorname{loc}}(\mathbb R^2)}\longrightarrow0,
\]
together with the uniform bound of $u_\varepsilon$ in
$L^\infty(0,T^*;(H^s(\R^2))^2)$, imply that this term converges to zero.
Hence, we deduce
\[
u_\varepsilon\cdot g_\varepsilon
\to
u\cdot g \text{ in } D'((0,T^*)\times\mathbb R^2).
\]
Finally, every weak-* accumulation point of the products equals $u\cdot g$, so the products converge along the whole sequence for which $u_\varepsilon$ has the stated weak-* limit.
\end{proof}

\begin{lemma}\label{lemma-0limit} It holds  $\omega_\varepsilon
\operatorname{div}u_\varepsilon  
\buildrel\ast\over\rightharpoonup 0$ in  $L^\infty(0,T^*;H^{s-2}(\mathbb{R}^2))$.
\end{lemma}
\begin{proof}
From \eqref{est-der-1} recall that $(\omega_\varepsilon
\operatorname{div}u_\varepsilon )_\varepsilon$ is bounded in  $L^\infty(0,T^*;H^{s-2}(\mathbb R^2)).$ By  the Banach–Alaoglu
theorem, it converges weakly-* (up to a subsequence) to some limit $\lambda\in L^\infty(0,T^*;H^{s-2}(\mathbb R^2)).$\\

    We show that $\lambda=0.$ Indeed,  multiplying \eqref{vorticity} by $\omega_\varepsilon$
    and using the identities
\[
\omega_\varepsilon\partial_t\omega_\varepsilon
=\frac12\partial_t|\omega_\varepsilon|^2,
\qquad
\omega_\varepsilon
u_\varepsilon\cdot\nabla\omega_\varepsilon
=
\frac12
u_\varepsilon\cdot\nabla|\omega_\varepsilon|^2,
\]
we obtain
\[
\frac12\partial_t|\omega_\varepsilon|^2
+\frac12
u_\varepsilon\cdot\nabla|\omega_\varepsilon|^2
+|\omega_\varepsilon|^2\operatorname{div}u_\varepsilon
+\frac{1}{\varepsilon}
\omega_\varepsilon\operatorname{div}u_\varepsilon
=0.
\]
Since
$\operatorname{div}\!\left(
\frac12|\omega_\varepsilon|^2u_\varepsilon
\right)=\frac12
u_\varepsilon\cdot\nabla|\omega_\varepsilon|^2+\frac12
|\omega_\varepsilon|^2\operatorname{div}u_\varepsilon,$
we obtain the following form
\[
\partial_t\!\left(\frac12|\omega_\varepsilon|^2\right)
+\operatorname{div}\!\left(
\frac12|\omega_\varepsilon|^2u_\varepsilon
\right)
+\frac12|\omega_\varepsilon|^2\operatorname{div}u_\varepsilon
+\frac{1}{\varepsilon}
\omega_\varepsilon\operatorname{div}u_\varepsilon
=0.
\]
Now, let  
$
\varphi\in C_c^\infty((0,T^*)\times\mathbb R^2),
$
we have
\begin{align*}
&\int_0^{T^*}\int_{\mathbb R^2}
\omega_\varepsilon\operatorname{div}u_\varepsilon
\varphi\,dxdt
\\&=\varepsilon\left[
\int_0^{T^*}\int_{\mathbb R^2}
\frac12|\omega_\varepsilon|^2\partial_t\varphi\,dxdt
+\int_0^{T^*}\int_{\mathbb R^2}
\frac12|\omega_\varepsilon|^2u_\varepsilon\cdot\nabla\varphi\,dxdt
-\int_0^{T^*}\int_{\mathbb R^2}
\frac12|\omega_\varepsilon|^2\operatorname{div}u_\varepsilon
\varphi\,dxdt
\right].
\end{align*}
By using that   $(\omega_\varepsilon)_\varepsilon$ is bounded in  $C(0,T^*;H^{s-1}(\mathbb R^2))$, $(u_\varepsilon)_\varepsilon$ is bounded in $C(0,T^*;(H^s(\mathbb R^2))^2)$ and $H^s(\mathbb R^2) \hookrightarrow W^{1,\infty}(\R^2)$, we deduce 
$$\int_0^{T^*}\int_{\mathbb R^2}
\omega_\varepsilon\operatorname{div}u_\varepsilon
\varphi\,dxdt \to 0 \text{ as } \varepsilon\downarrow 0,$$
and therefore $\lambda=0.$ Again, the uniqueness of the limit ensures the convergence of the whole sequence.
\end{proof}

\subsection{Step 4:  derivation of Hasegawa-Mima equation}
First, we present the following lemmas.
\begin{lemma}\label{Lem-ident-q}
     Let $0<\varepsilon\leq\varepsilon_0$ and $(\phi_\varepsilon,u_\varepsilon)$  
be a  family of solution to \eqref{Euler-rho-phi-Main}, given by \autoref{TH1-Main-BIS}, associated with the initial data $(\phi_0,u_0).$
Let $\phi \in  L^\infty(0,T^*;H^{s+1}(\mathbb R^2))$ such that 
\begin{align}\label{weak-limi-cv_*_2}
    \phi_\varepsilon\buildrel\ast\over\rightharpoonup \phi,
\text{ and }
u_\varepsilon \buildrel\ast\over\rightharpoonup -\nabla^\perp \phi 
\quad\text{in }
L^\infty(0,T^*;H^s(\mathbb R^2))  \text{ and }L^\infty(0,T^*;(H^s(\mathbb R^2))^2) \text{ resp.}
\end{align}
Then $q_\varepsilon=\phi_\varepsilon+ \nabla^\perp\cdot u_{\varepsilon}$, the unique solution to  \eqref{equ-derive-HM}, converges to $q=\phi-\Delta\phi$ in the following sense:
\begin{align*}
   q_\varepsilon &\buildrel\ast\over\rightharpoonup \phi-\Delta\phi
\quad\text{in }
L^\infty(0,T^*;H^{s-1}(\mathbb{R}^2)),\\
\notag   \partial_tq_\varepsilon &\buildrel\ast\over\rightharpoonup \partial_t(\phi-\Delta\phi)
\quad\text{in }
L^\infty(0,T^*;H^{s-2}(\mathbb{R}^2)).
\end{align*}
\end{lemma}
\begin{proof}
From \autoref{Lem-bound-q-appro}, $(q_\varepsilon)_\varepsilon$ is uniformly bounded in  $C(0,T^*;H^{s-1}(\mathbb{R}^2)) \cap C^1(0,T^*;H^{s-2}(\mathbb{R}^2))$. Therefore, by  the Banach–Alaoglu
theorem,   there exists $$q\in  L^\infty(0,T^*;H^{s-1}(\mathbb{R}^2)) \cap W^{1,\infty}(0,T^*;H^{s-2}(\mathbb{R}^2))$$  such that, up to a subsequence denoted by the same way, we get
\begin{align*}
   q_\varepsilon &\buildrel\ast\over\rightharpoonup q
\quad\text{in }
L^\infty(0,T^*;H^{s-1}(\mathbb{R}^2)),\\
   \partial_tq_\varepsilon &\buildrel\ast\over\rightharpoonup \partial_tq
\quad\text{in }
L^\infty(0,T^*;H^{s-2}(\mathbb{R}^2)).
\end{align*}It is worth remembering that, \textit{a priori}, there exists  $\kappa\in L^\infty(0,T^*;H^{s-2}(\mathbb{R}^2)) $ such that $\partial_tq_\varepsilon \buildrel\ast\over\rightharpoonup \kappa 
$ in $L^\infty(0,T^*;H^{s-2}(\mathbb{R}^2))$ and since the time derivative is understood in the sense of the distribution, we can identify the limit $\kappa=\partial_tq$ \textit{via} a  standard argument (see below for similar argument). Moreover, we have $q=\phi-\Delta\phi.$ Indeed, let $\varphi\in C_c^\infty((0,T^*)\times\mathbb{R}^2).$
Since
$q_\varepsilon=\phi_\varepsilon+ \nabla^\perp\cdot u_\varepsilon,$
we have
\[
\iint_{(0,T^*)\times\mathbb{R}^2}
q_\varepsilon\varphi\,dxdt
=
\iint_{(0,T^*)\times\mathbb{R}^2}
\phi_\varepsilon\varphi\,dxdt
+
\iint_{(0,T^*)\times\mathbb{R}^2}
(\nabla^\perp\cdot u_\varepsilon)\varphi\,dxdt.
\]
Using the identity
$-\displaystyle\int_{\mathbb{R}^2}(\nabla^\perp\cdot u)\,\varphi\,dx
=\int_{\mathbb{R}^2}
u\cdot\nabla^\perp \varphi \,dx,$
we obtain
\[
\iint_{(0,T^*)\times\mathbb{R}^2}
q_\varepsilon\varphi\,dxdt
=
\iint_{(0,T^*)\times\mathbb{R}^2}
\phi_\varepsilon\varphi\,dxdt
-
\iint_{(0,T^*)\times\mathbb{R}^2}
u_\varepsilon\cdot(\nabla\varphi)^\perp\,dxdt.
\]
Passing to the limit as $\varepsilon\to 0$ in the last identity and using the weak-$*$ convergences
\begin{align*}
    \phi_\varepsilon\buildrel *\over\rightharpoonup\phi
\quad\text{in }L^\infty(0,T^*;H^{s}(\mathbb{R}^2)) \text{ and } u_\varepsilon
\buildrel *\over\rightharpoonup
-\nabla^\perp \phi 
\quad\text{in }
L^\infty(0,T^*;(H^s(\mathbb{R}^2))^2),
\end{align*}
yields
\[
\iint_{(0,T^*)\times\mathbb{R}^2}
q\varphi\,dxdt
=
\iint_{(0,T^*)\times\mathbb{R}^2}
\phi\varphi\,dxdt
+
\iint_{(0,T^*)\times\mathbb{R}^2}
\nabla^\perp \phi \cdot(\nabla\varphi)^\perp\,dxdt.
\]
Since
$\nabla^\perp \phi \cdot(\nabla\varphi)^\perp
=\nabla\phi\cdot\nabla\varphi,
$ integrating by parts the second term once more, we get
\[
\iint_{(0,T^*)\times\mathbb{R}^2}
q\varphi\,dxdt
=
\iint_{(0,T^*)\times\mathbb{R}^2}
\left(
\phi-\Delta\phi
\right)\varphi\,dxdt.
\]
Since this identity holds for every
\(
\varphi\in C_c^\infty((0,T^*)\times\mathbb{R}^2),
\)
we conclude that $q=\phi-\Delta\phi.$
Due to the uniqueness of the limit, the convergence holds for the whole sequence, namely we proved 
\begin{align}\label{weak.cv-q}
   q_\varepsilon &\buildrel\ast\over\rightharpoonup \phi-\Delta\phi
\quad\text{in }
L^\infty(0,T^*;H^{s-1}(\mathbb{R}^2)),\\
\notag   \partial_tq_\varepsilon &\buildrel\ast\over\rightharpoonup \partial_t(\phi-\Delta\phi)
\quad\text{in }
L^\infty(0,T^*;H^{s-2}(\mathbb{R}^2)).
\end{align}
\end{proof}
Now, we have the following result about the strong convergence of $(q_\varepsilon)_\varepsilon.$ 
\begin{lemma}\label{strong-cv-lemma}
For every $\delta\in(0,1)$,
$q_\varepsilon\longrightarrow q$ strongly in $
C([0,T^*];H^{s-1-\delta}_{\mathrm{loc}}(\mathbb R^2)).$ In particular,
\begin{align}\label{initial-data-cv}\phi_0+ \nabla^\perp\cdot u_0 =(\phi-\Delta\phi)(0)
\end{align} and  $q_\varepsilon(t) \rightharpoonup q(t)$ in  $H^{s-1}(\mathbb{R}^2)$ for any $t\in [0,T^*[.$
\end{lemma}
\begin{proof}
Let $K\Subset\mathbb R^2$ be arbitrary, where $\Subset$ denotes the compact embedding. Choose
$\chi\in C_c^\infty(\mathbb R^2)$ such that
$\chi\equiv1$ on a neighborhood of $K$, and define $v_\varepsilon=\chi q_\varepsilon.$
Since multiplication by $\chi$ is continuous on Sobolev spaces, we have  $(v_\varepsilon)_\varepsilon$ is bounded in $
L^\infty(0,T^*;H^{s-1}(\mathbb R^2)),
$
and $\partial_tv_\varepsilon
=\chi\,\partial_tq_\varepsilon$
is bounded in
$L^\infty(0,T^*;H^{s-2}(\mathbb R^2)).$ Now,  let $D$ be a bounded  domain with   smooth boundary containing
$\operatorname{supp}\chi$. Then, one has
\begin{align*}
(v_\varepsilon)_\varepsilon
\text{ is bounded in }
L^\infty(0,T^*;H^{s-1}(D))
 \text{ and }
(\partial_tv_\varepsilon)_\varepsilon
\text{ is bounded in }
L^\infty(0,T^*;H^{s-2}(D)).
\end{align*}
 Let $\delta\in (0,1),$
since
$
H^{s-1}(D)
\Subset
H^{s-1-\delta}(D)
\hookrightarrow
H^{s-2}(D),
$
Simon's compactness theorem (see \textit{e.g.} \cite[Corollary 4]{Simoncompact}) implies that the family
$
(v_\varepsilon)_\varepsilon$ is relatively compact 
$
C([0,T^*];H^{s-1-\delta}(D)).
$
Consequently, after extraction of subsequence denoted by the same way,
 we have \[
v_\varepsilon
\longrightarrow
v
\quad\text{strongly in}\quad
C([0,T^*];H^{s-1-\delta}(D)) \text{ with } v\in L^\infty(0,T^*;H^{s-1}(D)).
\]

On the other hand, we have
$
q_\varepsilon
\buildrel *\over\rightharpoonup
q
\quad\text{in }
L^\infty(0,T^*;H^{s-1}(\mathbb R^2)),$
and multiplication by $\chi$ is continuous on
$H^{s-1}(\mathbb R^2)$. Hence, we get
$
v_\varepsilon=\chi q_\varepsilon
\buildrel *\over\rightharpoonup
\chi q
\quad\text{in }
L^\infty(0,T^*;H^{s-1}(\mathbb R^2)).
$
Therefore the strong limit satisfies
$
v=\chi q,
$ because of the uniqueness of the limit.
Thus, we have
\[
\chi q_\varepsilon
\longrightarrow
\chi q
\quad\text{strongly in}\quad
C([0,T^*];H^{s-1-\delta}(D)).
\]

Since $\chi\equiv1$ on $K$, we write
\[
q_\varepsilon
\longrightarrow
q
\quad\text{strongly in}\quad
C([0,T^*];H^{s-1-\delta}(K)).
\]
As $K\Subset\mathbb R^2$ was arbitrary, the following convergence holds
\begin{align}\label{strong-cv-time}
    q_\varepsilon
\longrightarrow
q
\quad\text{strongly in}\quad
C([0,T^*];H^{s-1-\delta}_{\mathrm{loc}}(\mathbb R^2)).
\end{align}

On the other hand, recall that for any $t\in [0,T^*[:$ $(q_\varepsilon(t))_\varepsilon$ is bounded in $H^{s-1}(\mathbb R^2))$ from \autoref{Lem-bound-q-appro}. Therefore, by  the Banach–Alaoglu theorem,   there exists $q_t\in H^{s-1}(\mathbb R^2))$ such that,
up to a subsequence denoted by the same way, we get $q_\varepsilon(t)\rightharpoonup q_t$ in $H^{s-1}(\mathbb R^2)).$ Finally, by using \eqref{strong-cv-time} one shows that $q_t=q(t)$ for any $t\in [0,T^*[.$
\end{proof}

Finally, we are in position to show that $\phi$ is the unique solution to Hasegawa-Mima equation, namely, to prove our main theorem.
\begin{proof}[Proof of \autoref{Theorem-main}]
    First, let $0<\varepsilon\leq\varepsilon_0$ and  we recall that $q_\varepsilon=\phi_\varepsilon+ \nabla^\perp\cdot u_{\varepsilon}$ is the unique solution to
    \begin{align*}
	\begin{cases}
    &\partial_t q_\varepsilon
+
u_\varepsilon\cdot\nabla q_\varepsilon
+ (\nabla^\perp\cdot u_{\varepsilon})
\operatorname{div}u_\varepsilon
+
u_\varepsilon\cdot g_\varepsilon
=0,\\[0.15cm]
&\big(\phi_\varepsilon, u_\varepsilon\big)    \big|_{t=0} = \big(\phi_0, u_{0}\big).
\end{cases}
\end{align*}
Let $\varphi\in C_c^\infty((0,T^*)\times\mathbb R^2).$
Testing the  above equation against $\varphi$ gives
\begin{align*}
&\iint_{(0,T^*)\times\mathbb{R}^2} \partial_tq_\varepsilon\varphi
 dxdt+\iint_{(0,T^*)\times\mathbb{R}^2}  u_\varepsilon\cdot\nabla q_\varepsilon\varphi dxdt
\\
&\quad
+\iint_{(0,T^*)\times\mathbb{R}^2}
(\nabla^\perp\cdot u_{\varepsilon})
\operatorname{div}u_\varepsilon\varphi dxdt
+
\iint_{(0,T^*)\times\mathbb{R}^2}
u_\varepsilon\cdot
 g_\varepsilon
\varphi dxdt=I^\varepsilon_1+I^\varepsilon_2+I^\varepsilon_3+I^\varepsilon_4 .
\end{align*}
Thanks to \autoref{Lem-ident-q}, we have 
\begin{align*}
 \lim_{\varepsilon\downarrow 0}I^\varepsilon_1= \lim_{\varepsilon\downarrow 0} \iint_{(0,T^*)\times\mathbb{R}^2} \partial_tq_\varepsilon\varphi
 dxdt=  \iint_{(0,T^*)\times\mathbb{R}^2} \partial_t(\phi-\Delta\phi)\varphi dxdt.
\end{align*}
Next, \autoref{lemma-0limit} ensures
\begin{align*}
 \lim_{\varepsilon\downarrow 0}I^\varepsilon_3= -\lim_{\varepsilon\downarrow 0} \iint_{(0,T^*)\times\mathbb{R}^2} ( \nabla^\perp\cdot u_{\varepsilon})
\operatorname{div}u_\varepsilon\varphi 
 dxdt=  0.
\end{align*}
Moreover, \autoref{Lem-log} gives
\begin{align*}
 \lim_{\varepsilon\downarrow 0}I^\varepsilon_4= \lim_{\varepsilon\downarrow 0} \iint_{(0,T^*)\times\mathbb{R}^2}
u_\varepsilon\cdot g_\varepsilon
\varphi dxdt&=\iint_{(0,T^*)\times\mathbb{R}^2}
u\cdot g
\varphi dxdt\\
&= -\iint_{(0,T^*)\times\mathbb{R}^2}
\nabla^\perp \phi \cdot g
\varphi dxdt.
\end{align*}
Finally, we have
\begin{align*}
 \lim_{\varepsilon\downarrow 0}I^\varepsilon_2= \lim_{\varepsilon\downarrow 0} \iint_{(0,T^*)\times\mathbb{R}^2}  u_\varepsilon\cdot\nabla q_\varepsilon\varphi dxdt=-\iint_{(0,T^*)\times\mathbb{R}^2}  \nabla^\perp \phi \cdot\nabla \big(\phi-\Delta\phi\big)\varphi dxdt.
\end{align*}
Indeed, 
since $u_\varepsilon\cdot\nabla q_\varepsilon
=
\operatorname{div}(u_\varepsilon q_\varepsilon)
-q_\varepsilon\operatorname{div}u_\varepsilon $ we have
\begin{align*}
     \lim_{\varepsilon\downarrow 0}I^\varepsilon_2&=  \lim_{\varepsilon\downarrow 0} \left[\iint_{(0,T^*)\times\mathbb{R}^2}  \operatorname{div}(u_\varepsilon q_\varepsilon)\varphi dxdt- \iint_{(0,T^*)\times\mathbb{R}^2}  q_\varepsilon\operatorname{div}u_\varepsilon\varphi dxdt\right]\\
     &=\lim_{\varepsilon\downarrow 0} \left[-\iint_{(0,T^*)\times\mathbb{R}^2}  q_\varepsilon u_\varepsilon\cdot\nabla\varphi dxdt- \iint_{(0,T^*)\times\mathbb{R}^2}  q_\varepsilon\operatorname{div}u_\varepsilon\varphi dxdt\right].
\end{align*}
Since $q_\varepsilon \to \phi-\Delta\phi$ strongly   in $C([0,T^*];H^{s-2}_{\mathrm{loc}}(\mathbb R^2))$, $u_\varepsilon \buildrel\ast\over\rightharpoonup -\nabla^\perp \phi $ in 
 $L^\infty(0,T^*;(H^s(\mathbb R^2))^2) $
  by \autoref{strong-cv-lemma} and   \autoref{lem-limits}, we obtain
  \begin{align*}
           \lim_{\varepsilon\downarrow 0}I^\varepsilon_2=\iint_{(0,T^*)\times\mathbb{R}^2}  (\phi-\Delta\phi) \nabla^\perp \phi \cdot\nabla\varphi dxdt,
  \end{align*}
  where we used the fact $\operatorname{div}(\nabla^\perp \phi )=0$ to deduce that
  \begin{align*}
      \lim_{\varepsilon\downarrow 0}  \iint_{(0,T^*)\times\mathbb{R}^2}  q_\varepsilon\operatorname{div}u_\varepsilon\varphi dxdt=0.
  \end{align*}
  Consequently, $\phi$ satisfies
  \begin{align*}
      \iint_{(0,T^*)\times\mathbb{R}^2} \partial_t(\phi-\Delta\phi)\varphi dxdt&+\iint_{(0,T^*)\times\mathbb{R}^2}  (\phi-\Delta\phi) \nabla^\perp \phi \cdot\nabla\varphi dxdt\\
      &-\iint_{(0,T^*)\times\mathbb{R}^2}
\nabla^\perp \phi \cdot g
\varphi dxdt=0 ,
  \end{align*}
  for any $\varphi\in C_c^\infty((0,T^*)\times\mathbb R^2).$ In other words, $\phi$ satisfies in the sense of distribution
  \begin{align}\label{limit-equ-dist}
      \partial_t(\phi-\Delta\phi)- \operatorname{div}\big[ (\phi-\Delta\phi) \nabla^\perp \phi \big]-
\nabla^\perp \phi \cdot g=0.
  \end{align}
  Since $\phi \in L^\infty(0,T^*;H^{s+1}(\mathbb R^2))$ and $  \partial_t(\phi-\Delta\phi)\in L^\infty(0,T^*;H^{s-2}(\mathbb{R}^2)) $  then \eqref{limit-equ-dist} holds in $L^\infty(0,T^*;H^{s-2}(\mathbb{R}^2)) $. Moreover, $\phi-\Delta\phi\in C_w(0,T^*;H^{s-1}(\mathbb{R}^2))$.
  In particular, $(\phi-\Delta\phi)(0) \in  H^{s-1}(\mathbb{R}^2)$ is well defined and due to the pointwise in time convergence of $q_\varepsilon$, we get $$(\phi-\Delta\phi)(0)=\phi_0+ \nabla^\perp\cdot u_0. $$
   In conclusion, $\phi$ solves the following Cauchy problem (Hasegawa-Mima equation)
    \begin{align}
	\begin{cases}
    & \partial_t(\phi-\Delta\phi)- \operatorname{div}\big[ (\phi-\Delta\phi) \nabla^\perp \phi \big]-
\nabla^\perp \phi \cdot g=0,\\[0.15cm]
&(\phi-\Delta\phi)\big|_{t=0} = \phi_0+ \nabla^\perp\cdot u_0,
\end{cases}
\end{align}
with $\phi \in L^\infty(0,T^*;H^{s+1}(\mathbb R^2))$, $\phi-\Delta\phi\in C_w(0,T^*;H^{s-1}(\mathbb{R}^2))$ and $$  \partial_t(\phi-\Delta\phi)\in L^\infty(0,T^*;H^{s-2}(\mathbb{R}^2)).$$
Since $\nabla^\perp \phi \cdot\nabla \phi=0$ and  $\operatorname{div}\nabla^\perp \phi =0$, we get
\begin{align}
	\begin{cases}
    & \partial_t(\phi-\Delta\phi)+\nabla^\perp \phi \cdot\nabla( \Delta\phi)-
\nabla^\perp \phi \cdot g=0,\\[0.15cm]
&(\phi-\Delta\phi)\big|_{t=0} = \phi_0+ \nabla^\perp\cdot u_0 \in H^{s-1}(\mathbb{R}^2).
\end{cases}
\end{align}
Since $\phi-\Delta\phi\in C_w(0,T^*;H^{s-1}(\mathbb{R}^2))$, by elliptic regularity theory  it follows that $$\phi\in C_w(0,T^*;H^{s+1}(\mathbb{R}^2)).$$
Finally, since $\phi\in  L^\infty(0,T^*;H^{s+1}(\mathbb R^2)) \cap C_w(0,T^*;H^{s+1}(\mathbb{R}^2))$ 
then $\phi$ is unique. Indeed, for two such solutions let $\psi=\phi_1-\phi_2$, $r=(I-\Delta)\psi$, $w=-\nabla^\perp\psi$ and $u_i=-\nabla^\perp\phi_i$. The difference equation is
\[
 \partial_t r+u_1\cdot\nabla r+w\cdot\nabla(\phi_2-\Delta\phi_2)+w\cdot g=0.
\]
Testing by $\psi$, the third term vanishes because $\operatorname{div}(\psi w)=0$. The drift term also vanishes:
\[
 \int_{\mathbb R^2}\psi w\cdot g\,dx
 =-\frac12\int_{\mathbb R^2}\nabla^\perp(\psi^2)\cdot\nabla G\,dx=0.
\]
Integration by parts in the remaining transport term gives
\[
 \frac{d}{dt}\big(\|\psi\|_2^2+\|\nabla\psi\|_2^2\big)
 \leq C\|\nabla u_1\|_\infty\big(\|\psi\|_2^2+\|\nabla\psi\|_2^2\big).
\]
The integrations can be justified by spatial cutoffs and approximation. Since $s>2$, Gronwall's inequality and the common initial data give $\psi=0$. Additionally,  due to the uniqueness of $\phi$ the convergence in \eqref{weak-limi-cv} holds for the whole sequence, which completes the proof.
\end{proof}
A direct corollary of the above proof is the following:
\begin{corollary}\label{cor-limit-HM}  
    There exists a unique $\phi\in  L^\infty(0,T^*;H^{s+1}(\mathbb R^2)) \cap C_w(0,T^*;H^{s+1}(\mathbb{R}^2))$ solution to the following Hasegawa-Mima equation \begin{align}\label{equ-lim-HM}
	\begin{cases}
    & \partial_t(\phi-\Delta\phi)+\nabla^\perp \phi \cdot\nabla( \Delta\phi)-
\nabla^\perp \phi \cdot  g=0,\\[0.15cm]
&(\phi-\Delta\phi)\big|_{t=0} = \phi_0+ \nabla^\perp\cdot u_0 \in H^{s-1}(\mathbb{R}^2).
\end{cases}
\end{align}
\end{corollary}
\section{Inverse energy cascade and drift waves in Hasegawa-Mima}\label{sec:cascade}

The purpose of this section is twofold. 
On one side it emphasizes the relevance, for the Hasegawa-Mima equation, of two
quantities that we may call \textit{generalized energy} and \textit{generalized enstrophy},
quantities which may suggest or explain the two main tricks developed in
\autoref{Subsection-Hs-est} and \autoref{subsec:vortcity}.
On the other side, we use them here to show a (conditional) property of
inverse energy cascade of the Hasegawa-Mima which explains its importance in
plasma turbulence theory, hence motivating the effort to reduce the complexity
of certain plasma models like the initial density-velocity system (where
detection of turbulent properties is more difficult).\\
If $\phi\in H^{1}\left(  \mathbb{R}^{2}\right)  $, the \textit{generalized
energy} $
E\left[  \phi\right]  =\displaystyle\int_{\mathbb{R}^{2}}\left(  \left\vert
\nabla\phi\right\vert ^{2}+\phi^{2}\right)  dx$
is well defined.

\begin{lemma}\label{Lem:conservation-energy} Let $r>3$ and 
 $\phi\in  L^\infty(0,T^*;H^{r}(\mathbb R^2)) \cap C_w(0,T^*;H^{r}(\mathbb{R}^2))$ be the unique solution to \eqref{equ-lim-HM}, given by \autoref{Theorem-main}. Then
 \begin{itemize}
     \item[(i)] we have $
E\left[  \phi\right]  =-\displaystyle\int_{\mathbb{R}^{2}}\left(  \Delta
\phi-\phi\right)  \phi dx.
$
\item[(ii)] Let $t\in (0,T^*),$  then $E\left[  \phi\left(  t\right)  \right]  $ is
constant \textit{i.e.}
$E\left[  \phi\left(  t\right)  \right]  =E\left[  \phi\left(  0\right)
\right].$
 \end{itemize}

\end{lemma}

\begin{proof}
(i) Since $\phi\in H^{r}\left(  \mathbb{R}^{2}\right)  $, all
terms in the identity above are finite and the identity follows by integration
by parts.

(ii) Since $\phi$ is sufficiently smooth for usual manipulations of calculus. Then%
\begin{align*}
\frac{d}{dt}\int_{\mathbb{R}^{2}}\left(  \left\vert \nabla
\phi\right\vert ^{2}+\phi^{2}\right)  dx &  =\int_{\mathbb{R}^{2}%
}\left(  2\nabla\phi\cdot\partial_{t}\nabla\phi
+2\phi\partial_{t}\phi\right)  dx  =-2\int_{\mathbb{R}^{2}}\phi\partial_{t}\left(  %
\Delta\phi-\phi\right)  dx\\
&  =-2\int_{\mathbb{R}^{2}}\phi \left[\nabla^\perp \phi \cdot\nabla( \Delta\phi)-
\nabla^\perp \phi \cdot
(\nabla\log n_0)  \right]  dx\\
&  =2\int_{\mathbb{R}^{2}}  \Delta\phi  (\nabla\phi)^{\perp}\cdot\nabla\phi dx+2\int_{\mathbb{R}^{2}}\phi 
\nabla^\perp \phi \cdot
(\nabla\log n_0) dx.
\end{align*}
Since $(\nabla\phi)^{\perp}\cdot\nabla\phi=0,$ one has $\displaystyle\int_{\mathbb{R}^{2}}  \Delta\phi  (\nabla\phi)^{\perp}\cdot\nabla\phi dx=0.$ On the other hand, we have
\begin{align*}
    \int_{\mathbb{R}^{2}}\phi 
\nabla^\perp \phi \cdot
(\nabla\log n_0)(0) dx&=-\int_{\mathbb{R}^{2}} 
(\nabla\phi^2)\cdot
(\nabla\log n_0)^\perp dx\\
&=\int_{\mathbb{R}^{2}} 
\phi^2
\operatorname{div}[(\nabla\log n_0)^\perp] dx=0.
\end{align*}
Therefore, we get 
$\displaystyle \frac{d}{dt}\int_{\mathbb{R}^{2}}\left(  \left\vert \nabla
\phi\right\vert ^{2}+\phi^{2}\right)  dx=0.
$
\end{proof}

If $\phi\in H^{2}\left(  \mathbb{R}^{2}\right)  $, then also the following
\textit{generalized enstrophy} is well defined:%
\[
\Omega\left[  \phi\right]  =\int_{\mathbb{R}^{2}}\left(  \left\vert \Delta\phi\right\vert ^{2}+\left\vert \nabla
\phi\right\vert ^{2}\right)  dx.
\]

\begin{remark}
Another version of it is
$\displaystyle\widetilde{\Omega}\left[  \phi\right]  =\int_{\mathbb{R}^{2}}\left( \Delta\phi-\phi\right)  ^{2}dx$
which make more clear the link with \autoref{subsec:vortcity}.
\end{remark}

 Let us recall the following.
\subsection{Generalized continuous Fourier basis}

Unlike the periodic case, where functions admit an expansion in terms of the
discrete Fourier modes $\{e^{ik\cdot x}\}_{k\in\mathbb Z^2}$, functions on
$\mathbb{R}^2$ are represented through the Fourier transform, which may be
viewed as a continuous analogue of the Fourier series.\\

For $\phi\in L^1(\mathbb{R}^2)$, the Fourier transform is defined by
\begin{equation}
\widehat{\phi}(\xi)
=
\frac{1}{2\pi}
\int_{\mathbb{R}^2}
\phi(x)e^{-ix\cdot\xi}\,dx,
\qquad
\xi\in\mathbb{R}^2,
\end{equation}
where $x\cdot\xi=x_1\xi_1+x_2\xi_2$ and  the inverse Fourier transform is given by
 $\phi(x)=\displaystyle\frac{1}{2\pi}
\int_{\mathbb{R}^2}
\widehat{\phi}(\xi)e^{ix\cdot\xi}\,d\xi.$
Plancherel theorem asserts that  Fourier transform extends uniquely to a unitary operator
\[
\mathcal{F}:L^2(\mathbb{R}^2)\rightarrow L^2(\mathbb{R}^2).
\] Consequently, the Fourier transform is
an isometry on $L^2(\mathbb{R}^2)$, namely
$\|\phi\|_{L^2(\mathbb{R}^2)}=\|\widehat{\phi}\|_{L^2(\mathbb{R}^2)},$
and  equivalently
\begin{equation}
(\phi,\psi)_{L^2(\mathbb{R}^2)}
=
(\widehat{\phi},\widehat{\psi})_{L^2(\mathbb{R}^2)},
\qquad
\forall\,\phi,\psi\in L^2(\mathbb{R}^2).
\end{equation}

Although the exponential functions
$e_\xi(x)=e^{ix\cdot\xi},$ $ \xi\in\mathbb{R}^2,$
do not belong to $L^2(\mathbb{R}^2)$ and therefore do not form an orthonormal
basis, they constitute a \textit{``continuous family of generalized Fourier mode''}.
Accordingly, every function $\phi\in L^2(\mathbb{R}^2)$ admits the
representation
\begin{equation}
\phi(x)
=
\frac{1}{2\pi}
\int_{\mathbb{R}^2}
\widehat{\phi}(\xi)e^{ix\cdot\xi}\,d\xi:=\int_{\R^2}\phi_k e_{k}\left(  x\right)dk,
\end{equation}
where the equality holds in the $L^2$-sense. Thus, the Fourier transform
provides a continuous decomposition of $L^2(\mathbb{R}^2)$, analogous to the
Fourier series expansion on the torus, with the discrete frequencies
$k\in\mathbb{Z}^2$ replaced by the continuous frequency variable
$\xi\in\mathbb{R}^2$.
Furthermore, differentiation is transformed into multiplication in Fourier
space: $
\widehat{\partial_j\phi}(\xi)
=
i\xi_j\widehat{\phi}(\xi),
\quad
j=1,2,$
and therefore
\begin{equation*}
\widehat{\nabla\phi}(\xi)
=
i\xi\,\widehat{\phi}(\xi),
\qquad
|\widehat{\nabla\phi}(\xi)|^2
=
|\xi|^2|\widehat{\phi}(\xi)|^2.
\end{equation*}
This property diagonalizes constant-coefficient differential operators.
In the sequel, we shall often use the \textit{generalized continuous Fourier basis} $e_{k}$,
$k\in\mathbb{R}^{2}$, of $L^{2}\left(  \mathbb{R}^{2}\right)  $. 

\subsection{On inverse cascade}

Using Fourier representation, we may write
\[
E\left[  \phi\right]  =\int_{\mathbb{R}^{2}}\left(  \left\vert
k\right\vert ^{2}+1\right)  \left\vert \phi_{k}\right\vert ^{2}dk
 \text{ and }
\Omega\left[  \phi\right]  =\int_{\mathbb{R}^{2}}\left(  \left\vert k\right\vert ^{4}+\left\vert k\right\vert
^{2}\right)  \left\vert \phi_{k}\right\vert ^{2}dk
.\] Introducing an energy
and enstrophy per Fourier component, called energy spectrum and enstrophy
spectrum, respectively,
\begin{align*}
E_{k}\left[  \phi\right]   &  =\left(  \left\vert k\right\vert
^{2}+1\right)  \left\vert \phi_{k}\right\vert ^{2} \text{ and }
\Omega_{k}\left[  \phi\right]     =\left(  \left\vert
k\right\vert ^{4}+\left\vert k\right\vert ^{2}\right)
\left\vert \phi_{k}\right\vert ^{2}.
\end{align*}
We have obtained the spectral decomposition of energy and enstrophy:
\[
E\left[  \phi\right]  =\int_{\mathbb{R}^{2}}E_{k}\left[  \phi\right]
dk,\qquad\Omega\left[  \phi\right]  =\int_{\mathbb{R}^{2}}\Omega_{k}\left[
\phi\right]  dk.
\]
Moreover, notice that
\[
\Omega_{k}\left[  \phi\right]  =E_{k}\left[  \phi\right]  \left\vert k\right\vert ^{2}.
\]

Now, consider a solution $\phi$ of the \eqref{equ-lim-HM} of class $H^{r}\left(
\mathbb{R}^{2}\right), r>3 $. Since $E\left[  \phi\left(  t\right)  \right]  $ is
constant, assuming $E\left[  \phi\left(  0\right)  \right]  >0$ (otherwise the
solution is trivial), we may introduce, for every $t\in\left[  0,T\right]  $,
the probability density%
\[
\varepsilon_{k}\left(  t\right)  :=\frac{E_{k}\left[  \phi\left(  t\right)
\right]  }{E\left[  \phi\left(  0\right)  \right]  }\qquad k\in\mathbb{R}%
^{2}\text{.}%
\]
It is the energy spectrum, normalized. Notice that, opposite to the energy
$E\left[  \phi\left(  t\right)  \right]  $ which is constant, the energy
spectrum and the probability density just introduced may change in time. 
\begin{align*}
    \text{The
crucial question we ask here is:\ how does it change?}
\end{align*}
 Assume for instance that
at time zero it is quite concentrated around some value of $k_{0}$; should we
expect it remain concentrated near $k_{0}$, or that it spreads, and in case it
spreads should we expect more shift towards larger or smaller values of $k$
w.r.t. $k_{0}$?

\begin{remark}
\label{remark spread in Fourier}We do not know how to prove that the spectrum
spreads, also because this is presumably true only often and in a statistical
sense. However, an indicator that some spread should occur is that the Fourier
formulation of \eqref{equ-lim-HM} (that we leave to the reader)\ contains a quadratic term
with infinitely many interactions. The nonlinear coupling has a tendency to
distribute ``information" among the Fourier components.
\end{remark}

One way to measure these properties is to introduce the following sort of
``average" and ``variance" of the energy spectrum

\[
m\left(  t\right)  :=\int_{\mathbb{R}^{2}}|k|\varepsilon_{k}\left(
t\right)  dk
\]%
\begin{align*}
\sigma^{2}\left(  t\right)   &  :=\int_{\mathbb{R}^{2}}\left\vert {|k|}%
-m\left(  t\right)  \right\vert ^{2}\varepsilon_{k}\left(  t\right)
dk=\int_{\mathbb{R}^{2}}\left\vert k\right\vert ^{2}%
\varepsilon_{k}\left(  t\right)  dk- m^2\left(  t\right)
\\
&  =\dfrac{\Omega\left[  \phi\left(  t\right)  \right]}{E(\phi(0))}  -m^2\left(  t\right).
\end{align*}

The main result of this section is:

\begin{theorem}\label{Thm:mean-var}
Let $\phi$ be the unique  solution of the \eqref{equ-lim-HM} of class $C_w\left(
\left[  0,T\right]  ;H^{r}\left(  \mathbb{R}^{2}\right)  \right)  $ with
$r>3$, as the one given by \autoref{Theorem-main}. Then $m\left(  t\right)  $ and
$\sigma^{2}\left(  t\right)  $ are differentiable and we have
\begin{align}\label{equ-mean-variance}
    \frac{d}{dt}m^2\left(  t\right)=-\frac{d}{dt}%
\sigma^{2}\left(  t\right)  -\dfrac{2}{E[\phi(0)]}\int_{\mathbb{R}^{2}}\left( \Delta\phi\right)  (\nabla\phi)^{\perp}
\cdot(\nabla\log n_{0})dx.
\end{align}

\end{theorem}

\begin{proof} Notice that $\phi$
belongs to  $C_w\left(  \left[  0,T\right]  ;H^{r}\left(
\mathbb{R}^{2}\right)  \right)  $ with $r>3$ and satisfies \eqref{equ-lim-HM}. This ensures that $\phi$ is  sufficiently smooth in
time  to make the expressions differentiable and
perform the computations in a classical sense. Therefore, we have
\begin{align*}
\frac{d}{dt}\Omega\left[  \phi\left(  t\right)  \right] &= 2\int
_{\mathbb{R}^{2}}\Delta\phi\partial_t\Delta\phi dx+2\int
_{\mathbb{R}^{2}}\nabla\phi\partial_t\nabla\phi dx\\ &  =2\int
_{\mathbb{R}^{2}}\left(  \Delta\phi\right)  \partial_{t}\left(
\Delta\phi-\phi\right)  dx\\
&  =2\int_{\mathbb{R}^{2}}   \Delta\phi  \left[  (
\nabla\phi)^{\perp}  \cdot\nabla\left(  \Delta\phi\right)  \right]  dx-2\int_{\mathbb{R}^{2}}   \Delta\phi   (
\nabla\phi)^{\perp}  \cdot(\nabla\log
n_{0})  dx.
\end{align*}
But%
\begin{align*}
\int_{\mathbb{R}^{2}}\left(  \Delta\phi\right)  \left[ 
(\nabla\phi)^{\perp}  \cdot\nabla \Delta\phi  \right]
dx &  =\frac{1}{2}\int_{\mathbb{R}^{2}}\left[   (\nabla
\phi)^{\perp}  \cdot\nabla\left(  \Delta\phi\right)  ^{2}\right]  dx\\
&  =-\frac{1}{2}\int_{\mathbb{R}^{2}} \operatorname{div}\left[ \left(
\nabla\phi\right)^{\perp}  \left(  \Delta\phi\right)  ^{2}\right]  dx=0.
\end{align*}
The identity  \eqref{equ-mean-variance} is proved, after using $\Omega(\phi(t))= E[\phi(0)]\big(\sigma^2(t)+m^2\left(  t\right)\big)$.
\end{proof}

The identity \eqref{equ-mean-variance} of \autoref{Thm:mean-var} is very rich in information. If the term
\[
-\int_{\mathbb{R}^{2}}\left(  \Delta\phi\right)  \left(
\nabla^{\perp}\phi\right)  \cdot(\nabla\log n_{0})dx
\]

is small and $\dfrac{d}{dt}\sigma^{2}\left(  t\right)  $ is large, namely
\[
\frac{d}{dt}\sigma^{2}\left(  t\right)  >-\dfrac{2}{E[\phi(0)]}\int_{\mathbb{R}^{2}}\left(\Delta\phi\right)  (\nabla\phi)^{\perp}
\cdot(\nabla\log n_{0})dx\]
(recall the heuristic discussion of \autoref{remark spread in Fourier} in
favor of the positivity of $\frac{d}{dt}\sigma^{2}\left(  t\right)  $), then
$m\left(  t\right) $ decreases. The decrease of
$ m\left(  t\right) $ corresponds to some form of shift
of the spectrum in the direction of small $k$. This is an indication of
\textit{energy transfer to large scales} (\textit{inverse cascade}, in loose
terms). And such energy transfer, seen in the physical space instead of
Fourier space, should correspond to the creation of larger energetic
structures, starting from smaller ones;\ the aggregation of small-scale
structures in larger ones.\\

A more speculative remark is that $\displaystyle\int_{\mathbb{R}^{2}}\left(  \Delta
\phi\right)  \left(  \nabla^{\perp}\phi\right)  \cdot(\nabla\log n_{0})dx$
is small also when $  (\nabla\phi)^{\perp}  \cdot(\nabla\log n_{0})$
is small and this occurs when zonal flows (elongated structures of $\phi$
perpendicular to $\nabla\log n_{0}$) are created. Thus under such
circumstances there is again a transfer of energy to large scales; and indeed,
zonal flows are maintained by small-scale energy transfer.
\subsection{Drift waves and \eqref{equ-lim-HM}}

This topic is very classical, see for instance \cite{HasMim,HorHas}.
We briefly recall it to see the power of  \eqref{equ-lim-HM} in comparison with the
continuity-momentum system we started from. The existence of drift waves (as
the waves described below are called)\ is again a central topic in the
understanding of plasma turbulence. \\

Waves are usually found by looking for complex-valued solutions of the form:
\begin{align*}
    \phi\left(  x,t\right)  =e^{i\left(  k\cdot x+\omega t\right)  } \text{ where } k\in\mathbb{R}^{2}, \omega\in\mathbb{R}.
\end{align*}
 To avoid complex-valued
expressions which are not coherent with the paper so far, let us look for
solutions of the form $\phi\left(  x,t\right)  =\sin\left(  k\cdot x+\omega
t\right)  $ (the result is analogous for $\phi\left(  x,t\right)  =\cos\left(
k\cdot x+\omega t\right)  $) and assume that $\nabla\log n_{0}$ is constant. We get the identity%

\[
\left(  \omega\left(  1+{\left\vert k\right\vert ^{2}}
\right)
+k^{\perp}\cdot\nabla\log n_{0}  \right)
\cos\left(  k\cdot x+\omega t\right)  =0.
\]
Therefore, if a pair $\left(  k,\omega\right)  $ satisfies the so-called
dispersion relation%
\[
\omega=-\frac{k  ^{\perp}\cdot\nabla\log n_{0}  }{1+\left\vert k\right\vert ^{2}}%
\]
then $\phi\left(  x,t\right)  =\sin\left(  k\cdot x+\omega t\right)  $ (and
also $\phi\left(  x,t\right)  =\cos\left(  k\cdot x+\omega t\right)  $) is a
solution; and there are no other pairs with this property.

\subsection*{Acknowledgements} 
We thank Emanuele Tassi and Philip J. Morrison, also authors of fundamental
works on the topic (see for instance \cite{TasChaMor}), for very useful
discussions and the indication of the importance of certain ordering
conditions. 
The research of F.F. and Y.T. is
funded by the European Union (ERC, NoisyFluid, No. 101053472). Views and
opinions expressed are however those of the authors only and do not necessarily
reflect those of the European Union or the European Research Council. Neither
the European Union nor the granting authority can be held responsible for them.

     \bibliography{HM-derivation}{}
\bibliographystyle{plain}

\end{document}